\documentclass[a4paper,12pt]{amsart}
\usepackage{amssymb}
\usepackage[]{latexsym}

\usepackage{amsthm}

\usepackage{lscape}

\usepackage[latin1]{inputenc}
\usepackage{color}
\usepackage[all,dvips]{xy}

\usepackage[dvips]{graphicx}

\newtheorem{thm}{Theorem}
\newtheorem{cor}[thm]{Corollary}
\newtheorem{lem}[thm]{Lemma}

\newtheorem{prop}[thm]{Proposition}

\theoremstyle{remark}
\newtheorem{rem}[thm]{Remark}
\newtheorem{exmp}[thm]{Example}

\theoremstyle{definition}

\newenvironment{pf}{\par\noindent{\bf Proof.}\enspace\ignorespaces}{\qed\par\par}

\def\qed{\hfill $\Box$}

\newcommand{\jac}[1]{{\rm Jac}( {#1} ) }

\newcommand{\pmodd}[1]{\,({\rm mod\,}{#1}) }

\newcommand{\red}{\mbox{red}}
\newcommand{\ord}{\mbox{ord}}

\newcommand{\bQ}{{\mathbb{Q}}}
\newcommand{\bF}{{\mathbb{F}}}
\newcommand{\bR}{{\mathbb{R}}}
\newcommand{\bP}{{\mathbb{P}}}
\newcommand{\bZ}{{\mathbb{Z}}}
\newcommand{\bA}{{\mathbb{A}}}

\begin{document}

\title{Five squares in arithmetic progression over quadratic fields}

\author[Enrique Gonz\'alez-Jim\'enez]{Enrique Gonz\'alez-Jim\'enez}
\address{Departamento de Matem{\'a}ticas, Universidad Aut{\'o}noma de Madrid;  and Instituto de Ciencias Matem{\'a}ticas (ICMat), 28049 Madrid, Spain}
\email{enrique.gonzalez.jimenez@uam.es}
\urladdr{http://www.uam.es/enrique.gonzalez.jimenez}

\author[Xavier Xarles]{Xavier Xarles}
\address{Departament de Matem\`atiques\\Universitat Aut\`onoma de
Barcelona\\08193 Bellaterra, Barcelona, Catalonia, Spain}
\email{xarles@mat.uab.cat}

\thanks{The first author was supported in part by grants MTM 2009--07291 (Ministerio de Ciencia e Innovaci\'on, Spain) and CCG08-UAM/ESP-3906 (Universidad Aut\'onoma de Madrid, Comunidad de Madrid, Spain). The second author was partially supported by the grant MTM 2009--10359 (Ministerio de Ciencia e Innovaci\'on, Spain).}


\maketitle

\begin{abstract}
We provide several criteria to show over which quadratic number
fields $\bQ(\sqrt{D})$ there is a non-constant arithmetic
progression of five squares. This is carried out by translating
the problem to the determination of when some genus five curves
$C_D$ defined over $\bQ$ have rational points, and then by using a
Mordell-Weil sieve argument, in addition to others. Using an
elliptic curve Chabauty-like method, we prove that the only
non-constant arithmetic progression of five squares over
$\bQ(\sqrt{409})$, up to equivalence, is $7^2, 13^2, 17^2, 409,
23^2$. Furthermore, we provide an algorithm allowing to construct
all the non-constant arithmetic progressions of five squares over
all quadratic fields. Finally, we state several problems and
conjectures related to this problem.
\end{abstract}

\section{Introduction}

A well-known result by Fermat, proved by Euler in 1780, states that
 four squares over $\bQ$ in arithmetic progression do not exist.
Recently, the second author showed that there are no six squares
in arithmetic progression over a quadratic field (see \cite{Xa}).
As a by-product of his proof, one reaches the conclusion that five
squares in arithmetic progression over quadratic fields exist, but
are all obtained from arithmetic progressions defined over $\bQ$.
The aim of this paper is to study over which quadratic fields
there are such five-square sequences, similar to how the first
author and J. Steuding studied the four-square sequences in
\cite{GJS}.

However, there is a great difference between the four-square
problem and the five-square problem: if a field contains four
squares in arithmetic progression, then it probably contains
infinitely many (non-equivalent modulo squares). But any number
field may only contain a finite number of five squares in
arithmetic progression. The reason for this is that the moduli space
parametrizing those objects is a curve of genus 5 (see section \ref{type}),
and may therefore only contain a finite number of points over a fixed
number field by Faltings' Theorem.

On the other hand, one can easily prove (remark \ref{inf}, section
\ref{sec9}), that there are infinitely many arithmetic progressions such
that their first five terms are squares over a quadratic field.
The conclusion is that there are infinitely many quadratic fields with
five squares in arithmetic progression.

In this paper, we will attempt to persuade the reader that, even though there are infinitely many such fields, they are few. For example, we
will show that there are only two number fields $\bQ(\sqrt{D})$,
for $D$ a square-free integer, with $D< 10^{13}$ having five
squares in arithmetic progression: those with $D=409$ and
$D=4688329$ (see Corollary \ref{computations}). In order to obtain
this result, we will develop a method, related to the Mordell-Weil
sieve, to prove that certain curves have no rational points.

The outline of the paper is as follows: in section \ref{sec2}, we
provide another proof of a result in \cite{Xa}, essential for our
paper. This result states that any arithmetic progression such
that its first five terms are squares over a quadratic field is
defined over $\bQ$. Using this result, we will show in section
\ref{type} that a number field $\bQ(\sqrt{D})$ contains five
different squares in arithmetic progression, if and only if, some
curve $C_D$ defined over $\bQ$ has $\bQ$-rational points. Next, we
study a little bit of the geometry of these curves $C_D$. In the
following sections, we provide several criteria to show when
$C_D(\bQ)$ is empty: when it has no points at $\bR$ or at $\bQ_p$
in section \ref{sec_local}, when it has an elliptic quotient of
rank $0$ in section \ref{sec_rank}, and when it does not pass some
kind of Mordell-Weil sieve in section \ref{sec_MW}. Section
\ref{sec_409} is devoted to computing all the rational points for
$C_{409}$. This is carried out by modifying the elliptic curve
Chabauty method, developed by Bruin in \cite{Br,Br0}. The result
obtained is that there are only $16$ rational points, all coming
from the arithmetic progression $7^2, 13^2, 17^2, 409, 23^2$.
Finally, in the last section, we give some tables related to the
computations, some values of $D$ where we do have rational points
in $C_D$, and we state several problems and conjectures.

\

{\bf Acknowledgements.}  We would like to thank Gonzalo Tornaria for aiding us with some
computations concerning the Corollary \ref{computations}. The authors thank the referees for helpful comments and suggestions.

\section{The $5$ squares condition}\label{sec2}

Recall that $n+1$ elements of a progression $a_0,\dots,a_n$ on a
field $K$ are in arithmetic progression if $a$ and $r
\in K$ exist, such that $a_i=a+i\cdot r$ for any $i=0,\dots,n$. This is
equivalent, of course, to having $a_i-a_{i-1}=r$ fixed for any
$i=1,\dots,n$. Observe that, in order to study squares in
arithmetic progression, we can and will identify the arithmetic
progressions $\{a_i\}$ and $\{a'_i\}$ such that an
$\alpha \in K^*$ with $a'_i=\alpha^2\, a_i$ for any $i$ exists. Hence, if $a_0\ne 0$,
we can divide all $a_i$ by $a_0$, and the corresponding
common difference is then $q=a_1/a_0-1$.

Let $K/\bQ$ be a quadratic extension. The aim of this section is
to show that any non-constant arithmetic progression whose first
five terms are squares over $K$, is defined  over $\bQ$ modulo the previous identification. Another proof of this result may be found in
\cite{Xa}.

First, let us consider the case of four squares in arithmetic progression over $K$.

\begin{prop} Let $K/\bQ$ be a quadratic extension, and let $x_i\in
K$ for $i=0,\dots,3$ be four elements, not all zero, such that
$x_i^2-x_{i-1}^2=x_j^2-x_{j-1}^2\in K$ for all $i,j=1,2,3$. Then
$x_0\ne 0$; and if we denote by $q:=(x_1/x_0)^2-1$, then $q=0$ or
$$\frac {(3q+2)^2}{q^2} \in \bQ.$$
\end{prop}

\begin{pf}
Observe that the conditions on $x_0,x_1,x_2, x_3$ are equivalent to the equations
$$x_0^2-2x_1^2+x_2^2=0 \ , \ x_1^2-2x_2^2+x_3^2=0,$$
which determine a curve $C$ in $\bP^3$. Observe also that $q$ is
invariant after multiplying all the $x_i$ by a constant, so we can
work with the corresponding point $[x_0:x_1:x_2:x_3]\in \bP^3$.
Using the previous equations, one shows easily that $x_0$ can not be
zero.

Before continuing, allow us to explain the strategy of the proof.
Since there are no four squares in arithmetic progression over
$\bQ$, the genus one curve $C$ satisfies
$C(\bQ)=\{[1:\pm1:\pm1:\pm1]\}$. Suppose we have a non-constant
map $\psi:C\rightarrow E'$ defined over $\bQ$, where $E'$ is an
elliptic curve defined over $\bQ$, such that $\psi(P)=0$ for all
$P\in C(\bQ)$. Denote by $\sigma$ the only automorphism of order
two of $K$, so $Gal(K/\bQ)=\{\sigma, id\}$. Then, for any point
$P\in C(K)$, $\psi(P)\oplus \psi(\sigma(P))$ must be $0$, so
$\psi(\sigma(P))=\sigma(\psi(P))=\ominus \psi(P)$. We will choose
such an elliptic curve $E'$ such that a Weierstrass equation
satisfies that the $x$-coordinate of $\psi(P)$ is equal to
$(3q+2)^2/q^2$. Since the $x$-coordinate is invariant by the
$\ominus$-involution, we will obtain the result.

Multiplying the equations $x_i^2=x_0^2+iq$,  for $i=1,2,3$ we obtain
$$(x_1x_2x_3)^2=(x_0^2+q)(x_0^2+2q)(x_0^2+3q).$$
So, replacing $q$ by $(x-2)x_0^2/6$, and $x_1x_2x_3/x_0^2$ by $y/6$, we get
the elliptic curve $E$ given by the equation
$$y^2=x^3+5x^2+4x,$$
with a map given by $f(x_0,x_1,x_2,x_3)=(2 x_3^2/x_0^2,6x_1x_2x_3/x_0^3)$. This
map is in fact an unramified degree four covering, corresponding
to one of the descendents in the standard $2$-descent. It sends the
8 trivial points to the points $(2,\pm6)$, which are torsion and
of order 4. We need a map that sends some trivial point to the
zero, so we just take $\tau(P):=P\oplus (2,-6)$. The map
$\tau:E\to E$ (not a morphism of elliptic curves) has equations
$$\tau(x,y)=\left(\frac
{2(x^2+14x+6y+4)}{(x-2)^2},-\frac{6(6xy+x^3+16x^2+32x+12y+8)}{(x-2)^3}\right).$$
The trivial points then go to the 0 point and the point $(0,0)$.

Now consider the standard 2-isogeny $\mu:E\to E'$, where $E'$ is
given by the equation $y^2=x^3-10x^2+9x$, given by
$$\mu(x,y)=\left(\frac{y^2}{x^2},\frac{y(4-x^2)}{x^2}\right)$$
(see for example \cite{Sil}, example III.4.5.).

The composition $\mu\circ\tau\circ f$ is exactly the map $\psi$ we
wanted. By applying the formulae above we obtain that the
$x$-coordinate of $\mu(\tau(f(x_0,x_1,x_2,x_3)))$ is exactly equal
to $\frac {(3q+2)^2}{q^2}$.
\end{pf}

\

We apply this proposition to obtain the result on five squares in
arithmetic progression.

\begin{cor}\label{cor2}Let $K/\bQ$ be a quadratic extension, and let $x_i\in
K$ for $i=0,\dots,4$ be five elements, not all zero, such that
$x_i^2-x_{i-1}^2=x_j^2-x_{j-1}^2\in K$ for all $i,j=1,2,3,4$. Then
$x_0\ne 0$, and if we denote by $q:=(x_1/x_0)^2-1$, then $q\in
\bQ$. In particular, $$\frac{x_i^2}{x_0^2}=1+i q\in\bQ, i=1,\dots,4.$$
\end{cor}

\begin{pf} Suppose $q\ne 0$. By the proposition we have that
$t_q:=(3q+2)^2/q^2 \in \bQ$ and that, if we
denote by $q':=(x_2/x_1)^2-1$, the same is true for $q'$. But
$q'=q/(q+1)$, so the condition for $q'$ is equivalent to
$t'_q:=(5q+2)^2/q^2\in \bQ$. But
$t'_q-t_q=16+8/q$, so $q\in \bQ$.
\end{pf}

\section{A diophantine problem over $\bQ$}\label{type}

Let $D$ be a square-free integer. We will say that two sets $S_1,S_2$ of $\bQ(\sqrt{D})$ are square equivalents if there exists $\alpha\in\bQ(\sqrt{D})$, $\alpha\ne 0$, such that $S_2=\alpha^2 S_1$. Notice that the previous equivalence is natural when the sets are formed by squares. Then, corollary \ref{cor2} showed that any arithmetic progression of $5$ squares over
$\bQ(\sqrt{D})$ is square equivalent to one arithmetic progression defined over $\bQ$.

\begin{lem}\label{lemanuevo}
Let $D$ be a square-free integer. Then an arithmetic progression of five squares over $\bQ(\sqrt{D})$ is square equivalent to one of the form $x_i^2=d_i X_i^2$ where $d_i=1$ or $D$, $X_i\in\bZ$ and the greatest common divisor of $x^2_0,\dots,x^2_4$ is square-free. We say that the $5$-term arithmetic progression is of type $I=\{i\,:\, d_i=D\}\subset\{0,\cdots,4\}$.
\end{lem}

\begin{pf}
Let $z_0,\dots,z_4\in \bQ(\sqrt{D})$ such that $z^2_0,\dots,z^2_4$ form an arithmetic progression. By Corollary \ref{cor2}, it is square equivalent to $y_i^2=1+i\,r/s$, $i=0,\dots,4$ for some $r,s\in\bZ$. In particular, it is square equivalent to $s^2 y_i^2=s^2+i s r$ with $s^2, s r\in\bZ$. Now, let $d$ be the greatest integer such that $d^2$ divides the greatest common divisor of $s^2 y_0^2,\dots,s^2 y^2_4$. Then the arithmetic progression $z_i^2$ is square equivalent to $x_i^2=(s/d)^2 y_i^2$, where the greatest common divisor of $x^2_0,\dots,x^2_4$ is square-free and since $x_i^2\in \bZ$ and $x_i\in \bQ(\sqrt{D})$ we have that $x_i^2=d_i X_i^2$ where $d_i=1$ or $D$ and $X_i\in\bZ$.
\end{pf}

Notice that $7^2,13^2,17^2,409,23^2$ is a $5$-term arithmetic progression over $\bQ(\sqrt{409})$ of type $\{3\}$, since $d_3=409$.

We define another equivalence relation on the set of $5$-term arithmetic progressions over
$\bQ(\sqrt{D})$ as follows: we say that two arithmetic progressions over $\bQ(\sqrt{D})$, $x_0^2,\dots,x_4^2$ and $y_0^2,\dots,y_4^2$ are equivalent if there exists $r\in\bQ$ and $\alpha=r^2$ or $\alpha=D\,r^2$ such that $y_i^2=\alpha x_i^2$ or $y^2_{4-i}=x^2_i$ for $i=0,\dots,4$.

\begin{lem}
A non-constant arithmetic progression of five squares over a
quadratic field, up to equivalences, is of type $\{3\}$.
\end{lem}

\begin{pf} Notice that up to the equivalences defined above, there are only few types of non-constant arithmetic
progressions of $5$ squares over quadratic fields: namely $\{i\}$ for $i=2,3,4$ and $\{i,j\}$ for $i=0,1$ and $j=1,\dots,4$ with $i<j$.

Now, assume that we have a $5$-term arithmetic progression $x_n^2=a+n q$, $n=0,\dots,4$, over $\bQ(\sqrt{D})$ of type $\{i,j\}$, then by Lemma \ref{lemanuevo} $x_i^2=D
X_i^2$, $x_j^2=D X_j^2$ and $x_k^2=X_k^2$ if $k\neq i,j$, where
$X_n\in \bZ$, $n=0,\dots,4$.  Let $p>3$ be  a prime dividing $D$.
Since $(j-i) q = x_j^2-x_i^2= D(X_j^2-X_i^2)$, we have $p | q$,
and therefore $p|a$. Thus we get that $p$ divides $x_n^2$ for all
$n=0,\dots,4$.

Let us see that, in fact, $p^2|x_n^2$ for all $n=0,\dots,4$, to
obtain a contradiction (recall that $x_n$ are not in $\bZ$, so
this is not automatic). Observe that for any $k\in\{0,\dots,4\}$
with $k\neq i,j$, we have that $x_k^2=X_k^2$ with $X_k\in \bZ$,
hence $p$ divides $X_k$ and so $p^2$ divides $x_k^2$. But now,
considering $k,l\in\{0,\dots,4\}$ such that $k,l\neq i,j$ and
$l>k$, we obtain that $(l-k) r = x_l^2-x_k^2$, and hence $p^2 | q$,
and therefore $p^2|a$. Then we have proved that the type $\{i,j\}$
is not possible over $\bQ(\sqrt{D})$ for $|D|>6$ and $|D|=5$. The cases $D=-6,-3,-2,-1,2$ and $3$ are not possible since there are no non-constant arithmetic progressions of four squares over $\bQ(\sqrt{D})$ (cf. \cite{GJS}). The remaining case $D=6$ is not possible but for a different argument, since there are infinitely many non-constant arithmetic progressions of four squares over $\bQ(\sqrt{6})$ (cf. \cite{GJS}).
We are going to prove that the types $\{i,j\}$ for $i=0,1$ and $j=1,\dots,4$ with $i<j$ are not possibles over $\bQ(\sqrt{6})$. Let define the following three conics in $\bP^2(\bQ)$:
$$
\begin{array}{l}
C_{1,i}\,:\,6X_i^2-12X_{i+1}^2+X_{i+2}^2=0,\\
C_{2,i}\,:\,6X_i^2-2X_{i+1}^2+6X_{i+2}^2=0,\\
C_{3,i}\,:\,X_i^2-2X_{i+1}^2+6X_{i+2}^2=0.
\end{array}
$$
Then it is straightforward to prove, using Hilbert symbols, that $C_{j,i}(\bQ)=\emptyset$ for $j=1,2,3$. Now, let be a $5$-term arithmetic progression $x_0^2,x_1^2,x_2^2,x_3^2,x_4^2$. Therefore $x_0,x_1,x_2,x_3,x_4$ are solutions of the system of equations
$$
x_0^2-2x_1^2+x_2^2=0,\quad x_1^2-2x_2^2+x_3^2=0,\quad x_2^2-2x_3^2+x_4^2=0.
$$
In particular, if this $5$-term arithmetic progression is over $\bQ(\sqrt{6})$ of type, say, $\{0,1\}$ then $x_0^2=6X_0^2$, $x_1^2=6X_1^2$ and $x_k^2=X_k^2$ for $k=2,3,4$ and $X_0,X_1,X_2,X_3,X_4\in\bQ$. Then the first equation of the previous system becomes $6X_0^2-12X_1^2+X_2^2=0$. That is, $[X_0:X_1:X_2]\in C_{1,0}(\bQ)$. But since $C_{1,0}(\bQ)=\emptyset$ we conclude that there are no non-constant $5$-term arithmetic progression over  $\bQ(\sqrt{6})$ of type $\{0,1\}$. For the remaining types we follow the same argument but replacing the conic $C_{1,0}$ with another conics as the following table shows: 
$$
\begin{array}{|c|c|c|c|c|c|c|}
\hline
\{0,1\} & \{0,2\} &\{0,3\} &\{0,4\} &\{1,2\} &\{1,3\} &\{1,4\}\\
\hline
C_{1,0} & C_{2,0} &C_{3,1} &C_{3,2} &C_{1,1} &C_{2,1} &C_{3,2}\\
\hline
\end{array}
$$

The type $\{4\}$ (or equivalently $\{0\}$) is not possible since
there are no non-constant arithmetic progressions of four squares
over the rationals.

To finish, let us see that the type $\{2\}$ is not possible. In
this case we have that $[x_0:x_1:x_3:x_4]\in \mathbb{P}^3(\bQ)$ is
a point on the intersection of two quadric surfaces in $\bP^3$:
$$
C_{\{2\}}\,:\,\left\{
\begin{array}{l}
X_1^2+2X_4^2-3X_3^2=0\\
X_3^2+2 X_0^2-3X_1^2=0.
\end{array}
\right.
$$

Note that the eight points $[ 1:\pm 1:\pm 1:\pm 1 ]$ belong to
$C_{\{2\}}$.  In the generic case the intersection of two quadric
surfaces in $\mathbb{P}^3$ gives an elliptic curve and, indeed,
this will turn out to be true in our case. A Weierstrass model for
this curve is given by $E:y^2=x(x+1)(x+9)$ (this is denoted by
\verb+48a3+ in Cremona's tables \cite{cremona,cremonaweb}). Using a computer
algebra package like \verb+MAGMA+ or \verb+SAGE+ (\cite{magma}, \cite{sage} resp.), we check that  $E(\bQ)\cong
\bZ/2\bZ\oplus\bZ/4\bZ $. Therefore  $C_{\{2\}}=\{[1:\pm 1:\pm 1:\pm 1 ]\}$, which implies $x_n^2=x_0^2$ for $n=0,1,3,4$. Deriving from that $Dx_2^2=x_0^2$ is then straightforward, which is impossible.
\end{pf}

\

Let $D$ be a square-free integer. We will denote by $C_D$ the
curve over $\bQ$ that classifies the arithmetic progressions of type
$\{3\}$. As a consequence of the previous result, we get the following geometric characterization.

\begin{cor}
Let $D$ be a square-free integer. Non-constant arithmetic
progressions of five squares over $\bQ(\sqrt{D})$, up to
equivalences, are in bijection with the set $C_D(\bQ)$.
\end{cor}

The curve $C_D$ has remarkable
properties that we are going to show in the sequel. First of
all, the curve $C_D$ is a non-singular curve over $\bQ$ of genus 5
that can be given by the following equations in $\bP^4$:
\begin{equation}\label{eqCD}
C_D\,:\,
\left\{
\begin{array}{l}
F_{012}:=X_0^2-2X_1^2+X_2^2=0\\
F_{123}:=X_1^2-2X_2^2+DX_3^2=0\\
F_{234}:=X_2^2-2DX_3^2+X_4^2=0
\end{array}\right.
\end{equation}
where we use the following convention: for $i,j,k
\in\{0,\dots,4\}$ distinct, $F_{ijk}$ denotes the curve that
classifies the arithmetic progressions $\{a_n\}_n$  (modulo
equivalence) such that $a_i=d_iX_i^2$,  $a_j=d_jX_j^2$,
$a_k=d_kX_k^2$, where $d_i=1$ if $i\ne 3$ and $d_3=D$.

Observe that we could describe also the curve $C_D$ by choosing
three equations $F_{ijk}$ with the only condition that all numbers
from $0$ to $4$ appear in the subindex of some $F_{ijk}$.

We have $5$ quotients of genus $1$ that are the intersection of
the two quadric surfaces in $\bP^3$ given by  $F_{ijk}=0$ and
$F_{ijl}=0$, where $i,j,k,l\in\{0,\dots,4\}$ distinct. Note that
these quotients consist on removing the variable $X_n$, where $n\ne
i,j,k,l$. We denote by $F^{(n)}_D$ to this genus $1$ curve.

These genus $1$ curves do not always have rational points (except for $F^{(3)}:=F^{(3)}_D$).  A Weierstrass model of the Jacobians of these genus $1$ curves can be computed by finding them in the case $D=1$
(using that $F^{(i)}_{1}$ has always some easy rational point) and
then twisting by $D$. Using the labeling of the Cremona's tables
\cite{cremona,cremonaweb}, one can check that $\jac{F^{(0)}_D}$ (resp.
$\jac{F^{(1)}_D}$, $\jac{F^{(2)}_D}$, $\jac{F^{(4)}_D}$) is the
$D$-twist of \verb+24a1+ (resp. \verb+192a2+, \verb+48a3+,
\verb+24a1+) and $\jac{F^{(3)}}$ is \verb+192a2+.  We denote by
$E^{(0)}$ (resp. $E^{(1)}$, $E^{(2)}$) the elliptic curve
\verb+24a1+ (resp. \verb+192a2+, \verb+48a3+) and by $E_D^{(i)}$
the $D$-twist of $E^{(i)}$, for $i=0,1,2$. Observe
also that $E^{(2)}=E^{(0)}_{-1}$, so $E^{(2)}_D=E^{(0)}_{-D}$.

Note that, in particular, we have shown the following result about
the decomposition in the $\bQ$-isogeny class of the Jacobian of
$C_D$.

\begin{lem}
Let $D$ be a square-free integer. Then
$$
\jac{C_D}\stackrel{\bQ}{\sim} \left(E_D^{(0)}\right)^2 \times
E_D^{(2)} \times E_D^{(1)} \times E^{(1)}\,.
$$
\end{lem}

\section{Local solubility for the curve $C_D$}\label{sec_local}

The aim of this section is to describe under which conditions with
respect to $D$ the curve $C_D$ has points in $\bR$ and $\bQ_p$ for
all prime numbers $p$.

\begin{prop}\label{local}
Let $D$ be a square-free integer. Then $C_D$ has points
in $\bR$ and in $\bQ_p$ for all primes $p$ if and only if $D>0$, 
$D\equiv \pm 1 \pmodd 5$ and for all primes $p$
dividing $D$, $p \equiv 1 \pmodd{24}$.
\end{prop}

This result is deduced from the following lemmas.

\begin{lem}
Let $D$ be a square-free integer. The curve $C_D$ has points in $K$, for $K=\bR$, $\bQ_2$, $\bQ_3$ and $\bQ_5$ if and only if $D$ is square in $K$. Explicitly, $D>0$,  $D\equiv 1\pmodd{8}$, $D\equiv
1 \pmodd{3}$ and $D\equiv \pm 1 \pmodd 5$, respectively.
\end{lem}

\begin{pf}
First, suppose that $D$ is a square over a field $K$. Then the curve $C_D$ contains the following sixteen points $[1:\pm 1: \pm 1: \pm 1/\sqrt{D}:\pm 1]$. This shows one of the implications. In order to show the other implication we will consider the different fields separately. Suppose that $C_D(K)\neq\emptyset$.

If $K=\bR$, the equation $F_{234}=0$ implies that $2DX_3^2=X_2^2+X_4^2$, which has solutions in $K$ only if $D>0$.

Consider now the case $K=\bQ_2$. On one hand, the conic given by
the equation $F_{123}=X_1^2-2X_2^2+DX_3^2$ has points in $\bQ_2$
if and only if $(2,-D)_2=1$, where $(\,,\,)_2$ denotes the Hilbert
symbol. This last condition is equivalent to $D\equiv \pm 1 \pmodd
8$ or $D\equiv \pm 2 \pmodd {16}$. On the other hand, doing the
same argument for the equation $F_{234}=X_2^2-2DX_3^2+X_4^2$ we
get the condition $(-1,2D)_2=1$, which implies $D\equiv 1 \pmodd
4$ or $D\equiv 2 \pmodd 8$. So we get $D$ odd and $D\equiv 1
\pmodd 8$, or $D$ even and $D\equiv 2 \pmodd {16}$. This last case
is equivalent, modulo squares, to the case $D=2$ and it is easy to
show that $C_2(\bQ_2)=\emptyset$.

If $K=\bQ_3$, considering the reduction modulo $3$ of the conic given by the equation $F_{023}=0$ we obtain that  $D\not\equiv -1 \pmodd 3$. Similarly, we have $D\not\equiv 0 \pmodd 3$ using $F_{123}=0$.

Finally if $K=\bQ_5$, one can show by an exhaustive search that there is no point in $C_D(\bF_5)$ if $D\equiv \pm 2 \pmodd 5$. The case $D\equiv 0 \pmodd 5$ is discharged by using $F_{123}=0$ modulo 5.
\end{pf}

In the following we will study the remaining primes $p>5$ in two
separate cases, depending if $p$ divides or not $D$. The first
observation is that the case that $p$ does not divide $D$
corresponds to the good reduction case.

\begin{lem} Let $p>3$ be a prime not dividing $D$. Then the model of $C_D$ given by the equations
$F_{012}$, $F_{123}$ and $F_{234}$ has good reduction at $p$.
\end{lem}

\begin{pf} We use the Jacobian criterion.
The Jacobian matrix of the system of equations defining $C_{D}$ is
$$A:=( \partial F_{i(i+1)(i+2)}(X_i,X_{i+1},X_{i+2})/\partial X_j)_{0\le i\le 2,0\le j\le 4}.$$

For any $j_1<j_2$, denote by $A_{j_1,j_2}$ the square matrices obtained
from $A$ by deleting the columns $j_1$ and $j_2$. Their determinants are equal to
$$|A_{j_1,j_2}|=k_{j_1,j_2}\cdot \prod_{i\ne j_1,j_2} X_i $$
where
$$k_{0,1}=2^3D\ , \ k_{0,2}=-2^4D \ , \ k_{0,3}=2^33 \ , \  k_{0,4}=2^5D \ , \  k_{1,2}=2^3D,$$
$$ k_{1,3}=-2^4 \ , \ k_{1,4}=2^3 3 D \ , \ k_{2,3}=2^3 \ , \
k_{2,4}=-2^4 D\ , \ k_{3,4}=2^3.$$

Now, suppose we have a singular point of $C_D(\bF_p)$. Then the
matrix $A$ must have rank less than $3$ evaluated at this point, so all
these determinants must be 0. But, if $p>3$ and does not divide
$D$, then all the products of three homogeneous coordinates must
be zero, so the point must have three coordinates equal to 0, which
is impossible again if $p>3$. \end{pf}

\begin{lem} Let $p>5$ be a prime such that $p$ does not divide $D$.
 Then $C_D(\bQ_p)\ne \emptyset$. \end{lem}

\begin{pf} First of all, by Hensel's lemma, and since $C_D$ has
good reduction at $p$, we have that any point modulo $p$ lifts to
some point in $\bQ_p$. So we only need to show that $C_D(\bF_p)\ne
\emptyset$. Now, because of the Weil bounds, we know that $\sharp
C_D(\bF_p)> p+1 - 10 \sqrt{p}$. So, if $p>97$, then $C_D(\bF_p)\ne
\emptyset$ and we are done. For the rest of primes $p$, $5<p<97$,
an exhaustive search proves the result.
\end{pf}

\

We suspect that there should be some reason, besides the Weil
bound, that for all primes $p>5$ not dividing $D$, the curve $C_D$
has points modulo p, that should be related to the special form it
has or to the moduli problem it classifies.

\begin{lem}\label{localp} Let $p$ be a prime dividing $D$, and $p>3$.
Then $C_D(\bQ_p)\ne \emptyset$ if and only if $p\equiv 1
\pmodd{24}$.
\end{lem}
\begin{pf} We will show that a necessary and sufficient condition for
$C_D(\bQ_p)\ne \emptyset$ is that $2$, $3$ and $-1$ are all
squares in $\bF_p$. This happens exactly when $p\equiv 1
\pmodd{24}$. Note that this condition is sufficient since  $[\sqrt{3}:\sqrt{2}:1: 0:\sqrt{-1}]$ belongs to $C_D$.

Suppose that we have a point in $C_D(\bQ_p)$ given by a solution
of the equations $F_{ijk}$ in projective coordinates
$[x_0:x_1:x_2:x_3:x_4]$, with $x_i\in\bZ_p$, and such that not all
of them are divisible by $p$. The first observation is that only
one of the $x_i$ may be divisible by $p$; since if two of them,
$x_i$ and $x_j$, are divisible by $p$, we can use the equations
$F_{ijk}$ in order to show that $x_k$ is also divisible, for any
$k$.

Now, reducing $F_{123}$ modulo $p$, we obtain that $2$ must be a
square modulo $p$. Reducing $F_{234}$ modulo $p$ we obtain that $-1$ must be a
square modulo $p$. And finally, reducing
$F_{034}=X_0^2-4DX_3^2+3X_4^2$ modulo $p$ we obtain that $3$ must be
a square modulo $p$. So the conditions are necessary.
\end{pf}

\section{The rank condition}\label{sec_rank}

Let us begin by recalling the well-known 2-descent on elliptic curves,
as explained for example in \cite[Chapter X, Prop. 1.4]{Sil}. Consider $E$ an elliptic curve over a number field $K$ given by an
equation of the form
$$y^2=x(x-e_1)(x-e_2)\ , \ \mbox{ with }e_1,e_2\in K.$$
Let $S$ be the set of all archimedean places, all places dividing
2 and all places where $E$ has bad reduction. Let $K(S,2)$ be the
set of all elements $b$ in $K^*/K^{*2}$ with $\ord_v(b)$ is even for all
$v\notin S$. Given any $(b_1,b_2)\in K(S,2)\times K(S,2)$,
define the curve $H_{b_1,b_2}$ given as intersection of two
quadrics in $\bP^3$ by the equations
$$
H_{b_1,b_2}:
\left\{
\begin{array}{l}
b_1z_1^2-b_2z_2^2=e_1z_0^2, \\
b_1z_1^2-b_1b_2z_3^2=e_2z_0^2.
\end{array}
\right.$$
Then the curves $H_{b_1,b_2}$ do not depend on the representatives, up to isomorphism, and they have genus one with Jacobian $E$. Moreover, we have a natural degree
four map $\phi_{b_1,b_2}: H_{b_1,b_2}\to E$ given by
$$\phi_{b_1,b_2}(z_0,z_1,z_2,z_3):=(b_1(z_1/z_0)^2,b_1b_2z_1z_2z_3/z_0^3).$$

Moreover, the 2-Selmer group $S^{(2)}(E/K)$  of $E$ may be
identified with the subset
$$
S^{(2)}(E/K)=\{(b_1,b_2)\in K(S,2)\times K(S,2)\, |\,
H_{b_1,b_2}(K_v)\ne \emptyset\  \forall v \mbox{ place in $K$}\}.
$$
The group $E(K)/2E(K)$ may be described, via the natural injective map  $\psi:E(K)/2E(K) \to S^{(2)}(E/K)$ defined by
$$
\psi((x,y))=
\left\{
\begin{array}{ccl}
(x,x-e_1) & \mbox{if} & x\ne 0,e_1\\
(e_2/e_1,-e_1) & \mbox{if} & (x,y)=(0,0)\\
(e_1,(e_1-e_2)/e_1) & \mbox{if} & (x,y)=(e_1,0)\\
\end{array}
\right.
$$
and $\psi(0)=(1,1)$, as the subgroup consisting of $(b_1,b_2)\in K(S,2)\times K(S,2)$ such that $H_{b_1,b_2}(K)\ne \emptyset$.

The following lemma is elementary by using the description above,
and it is left to the reader.

\begin{lem}\label{Selmer}
 Let $H$ be a genus 1 curve over a number field $K$ given by an equation of the form
$$
H:
\left\{
\begin{array}{l}
b_1z_1^2-b_2z_2^2=e_1z_0^2 \\
b_1z_1^2-b_1b_2z_3^2=e_2z_0^2
\end{array}
\right.$$
 for some $b_1,b_2,e_1,e_2 \in K$.
Let $D\in K^*$ and consider the curves $H_D^{(1)}$, $H_D^{(2)}$
and $H_D^{(3)}$ given by replacing $z_1^2$ by
$Dz_1^2$,  $z_2^2$ by $Dz_2^2$ and $z_3^2$ by $Dz_3^2$ respectively in the equations above. Then $H_D^{(1)}$, $H_D^{(2)}$
and $H_D^{(3)}$ are homogeneous spaces for  the elliptic curve
$E_D$, the twist by $D$ of $E$, given by the Weierstrass equation
$y^2=x(x-De_1)(x-De_2)$.

Moreover, if $S_D$ denotes the set of all archimedean places, all
places dividing $2D$ and all places where $E$ has bad reduction,
the curves $H_D^{(1)}$, $H_D^{(2)}$ and $H_D^{(3)}$ correspond
respectively to the elements $(Db_1,b_2)$, $(b_1,Db_2)$ and
$(Db_1,Db_2)$ in $K(S_D,2)\times K(S_D,2)$.
 \end{lem}

\begin{prop} \label{rank}
Let $D>0$ be a square-free integer. A necessary condition for the
existence of 5 non-trivial squares in arithmetic progression over
$\bQ(\sqrt{D})$ is that the elliptic curves $E_D^{(0)}$ and
$E_D^{(2)}$ given by equations $Dy^2=x(x+1)(x+4)$ and
$Dy^2=x(x+1)(x+9)$ have rank 2 or more over $\bQ$, and that the
elliptic curve $E_D^{(1)}$ given by the equation
$Dy^2=x(x+2)(x+6)$ has positive rank.
\end{prop}

\begin{pf} Assume we have 5 non-trivial squares
in arithmetic progression over $\bQ(\sqrt{D})$. By using the results from
section \ref{type}, we can assume that such squares are of the
form $x_0^2$, $x_1^2$, $x_2^2$, $Dx_3^2$ and $x_4^2$, with
$x_i\in\bZ$. The condition of being in arithmetic progression is
equivalent to $x_0^2=a$ $x_1^2=a+q$,  $x_2^2=a+2q$,  $Dx_3^2=a+3q$
and $x_4^2=a+4q$ for some $a, q\in \bZ$. From these equations we
easily obtain that the following homogeneous spaces attached to $E_D^{(0)}$ have rational points:
$$
\left\{
\begin{array}{l}
 2(DX_3)^2-3DX_2^2=-DX_0^2\\
 2(DX_3)^2-6DX_1^2=-4DX_0^2
\end{array}
\right.
 \quad\mbox{and}\quad
 \left\{
\begin{array}{l}
2DX_4^2-3(DX_3)^2=-DX_1^2\\
2DX_4^2-6DX_2^2=-4DX_1^2
\end{array}
\right.
$$
which give $(2,3D)$ and $(2D,3)\in S^{(2)}(E_D^{(0)}/\bQ)$ by using Lemma \ref{Selmer}. Since
we are supposing both curves have points in $\bQ$, they correspond
to two points: $P_1$ and $P_2$ in $E_D^{(0)}(\bQ)$. In order to
show they have infinite order, we only need to show that the
symbols $(2,3D)$ and $(2D,3)$ are not in
$$\psi(E_D^{(0)}[2])=\{(1,1),(4,4D)=(1,D),(-D,-1), (-D,-D)\}$$
which is clear since $D>0$. In order to show that $P_1$ and $P_2$
are independent modulo torsion, it is sufficient to show that
$(2,3D)(2D,3)=(D,D)$ is not in $\psi(E_D^{(0)}[2])$, which is
again clear. So $E_D^{(0)}(\bQ)$ has rank $>1$.

The other conditions appear in similar fashion. We have
 $$
 \left\{
\begin{array}{l}
 3DX_4^2-4(DX_3)^2=-DX_0^2\\
  3DX_4^2-12DX_1^2=-9DX_0^2
\end{array}
\right.
\quad\mbox{and}\quad
\left\{
\begin{array}{l}
3DX_0^2-4DX_1^2=-DX_4^2\\
3DX_0^2-12D^2X_3^2=-9DX_4^2
\end{array}
\right.
$$
 which give $(3D,1)$ and $(3D,4D)=(3D,D) \in S^{(2)}(E_D^{(2)}/\bQ)$,
 again giving two independent points in $E_D^{(2)}(\bQ)$.

 Finally, we have
 $$ 6DX_4^2-2(2DX_3)^2=-2DX_0^2\ , \  6DX_4^2-12DX_1^2=-6DX_0^2$$
 which gives $(6D,2) \in S^{(2)}(E_D^{(1)}/\bQ)$, giving a non-torsion point in $E_D^{(1)}(\bQ)$.
 \end{pf}

\begin{rem}\label{rank_rem} Suppose that $D$ verifies the conditions in Proposition
\ref{local}, so $C_D(\bQ_p)\ne \emptyset$ for all $p$. Then the
root number of $E_D^{(0)}$ and $E_D^{(2)}$ is $1$ independently of
$D$ in both cases, and the root number of $E_D^{(1)}$ is always
$-1$. This is because the root number of the twist by $D$ of
an elliptic curve $E$ of conductor $N$, if $N$ and $D\equiv 1
\pmod{4}$ are coprime, is equal to the Kronecker
symbol $(D/-N)$ times the root number of $E$ (see, for example, section 4.3 of \cite{SiC}, which is deduced from
the Corollary to Proposition 10 in \cite{Ro})
In our case, and assuming $D$ verifies the conditions
in Proposition \ref{local}, we obtain that the root number of $E_D^{(i)}$
is equal to the root number of $E^{(i)}$, since $(D/-N)=1$ for
$N=24,48,192$.

 Assuming the Parity conjecture, this implies that the rank of
$E_D^{(0)}$ and $E_D^{(2)}$ is always even, and the rank of
$E_D^{(1)}$ is always odd. So the last condition in the
proposition is (conjecturally) empty.
\end{rem}

\subsection{Ternary Quadratic Forms}

It has been shown in Proposition \ref{rank} that a necessary
condition for the existence of a non-constant arithmetic
progression of $5$ squares over a quadratic field $\bQ(\sqrt{D})$
is that the elliptic curve $E_D^{(0)}$ and $E_D^{(2)}$ have
ranks $\ge 2$. In this part, we would like to describe some
explicit results concerning the ranks of these curves, thereby obtaining
an explicit computable condition.

\begin{rem}
The elliptic curve $E_D^{(0)}$ (resp. $E_D^{(2)}$) parametrizes
non-constant arithmetic progression of $4$ squares over
$\bQ(\sqrt{D})$ (resp. $\bQ(\sqrt{-D})$) (cf. \cite{GJS}).
Therefore, a necessary condition for the existence of a non-constant
arithmetic progression of $5$ squares over $\bQ(\sqrt{D})$ is the
existence of a non-constant arithmetic progression of $4$ squares
over $\bQ(\sqrt{D})$ and over $\bQ(\sqrt{-D})$.
\end{rem}

Using Waldspurger's results and Shimura's correspondence {\it a
la} Tunnell, several results were obtained by Yoshida \cite{yoshida2} on
the ranks of $E_D^{(0)}$ and $E_D^{(2)}$. In particular, we use
his results corresponding to the case  $D\equiv 1 \pmodd {24}$ in order to
apply them to our problem.

\begin{prop}\label{ternary}
Let $D$ be a square-free integer. If $Q(x,y,z)\in\mathbb{Z}[x,y,z]$
is a ternary quadratic form, denote by $r(D,Q(x,y,z)) $ the
number of integer representations of $D$ by $Q$. If
$$
\begin{array}{rcl}
r(D,x^2+12y^2+15z^2+12yz) & \neq  & r(D,3x^2+4y^2+13z^2+4yz)\\
& \mbox{or} & \\
 r(D,x^2+3y^2+144z^2) & \neq & r(D,3x^2+9y^2+16z^2),
\end{array}
$$
then there are no non-constant arithmetic progressions of $5$ squares over $\bQ(\sqrt{D})$.
\end{prop}

\begin{pf}
First of all, by the Proposition \ref{local} we have that $D\equiv
1\pmodd{24}$. Now, Yoshida constructs two cusp forms
of weight $3/2$ denoted by $\Phi_{3,-3}$ and $\Phi_{1,1}$, such
that if we denote by $a_D(\Phi_{3,-3})$ (resp.  $a_D(\Phi_{1,1})$)
the $D$-th coefficient of the Fourier $q$-expansion of
$\Phi_{3,-3}$ (resp.  $\Phi_{1,1}$), we have
$$
\begin{array}{l}
 \mbox{$a_D(\Phi_{3,-3})=0$ if and only if $L(E^{(0)}_D,1)=0$,}\\[1mm]
 \mbox{$a_D(\Phi_{1,1})=0$ if and only if $L(E^{(2)}_D,1)=0$.}\\
\end{array}
$$
Then by the definition of these cusp forms we have:
$$
\begin{array}{ll}
 \!\!a_D(\Phi_{3,-3}) \!\!&\!\! = r(D,x^2+12y^2+15z^2+12yz)-r(D,3x^2+4y^2+13z^2+4yz)\,,\\
 \!\!a_D(\Phi_{1,1})  \!\!&\!\! =  r(D,x^2+3y^2+144z^2)-r(D,3x^2+9y^2+16z^2)\,,
\end{array}
$$
which concludes the proof.
\end{pf}

\begin{rem}For $D=2521$, the conditions in Propositions \ref{local}, \ref{rank} and \ref{ternary} are fulfilled, and in fact all of the
relevant genus $1$ curves have rational points. But we will show in Corollary \ref{computations} that $C_{2521}(\bQ)=\emptyset$.
 \end{rem}

\section{The Mordell-Weil sieve}\label{sec_MW}

In this section we want to develop a method to test when $C_D$ has
no rational points based on the Mordell-Weil sieve (see \cite{Scharaschkin}, \cite{Flynn}, \cite{Poonen}, \cite{Stoll}, \cite{Bruin:Stoll}).

The idea is the following: Suppose we have a curve $C$ defined
over a number field $K$ together with a map $\phi:C\to A$ to an abelian variety $A$ defined
over $K$. We want to show that $C(K)=\emptyset$, and we know that
$\phi(C(K))\subset H\subset A(K)$, a certain subset of $A(K)$. Let
$\wp$ be a prime of $K$ and consider the reduction at $\wp$ of all
the objects $\phi_{\wp}: C_{\wp} \to A_{\wp}$, together with the
reduction maps $\red_{\wp}: A(K) \to A(k_{\wp})$, where $k_{\wp}$
is the residue field at $\wp$. Now, we have that $\red_{\wp}(C(K))
\subset \phi_{\wp}(C(k_{\wp})) \cap \red_{\wp}(H)$, so
$$
\phi(C(K))\subset
H^{(\wp)}:=\red_{\wp}^{-1}\Big(\phi_{\wp}(C(k_{\wp})) \cap
\red_{\wp}(H)\Big).$$ By considering enough primes,
it could occur that $$\bigcap_{\mbox{some primes }\wp}
\!\!\!\!\!\!\!\!\!H^{(\wp)}=\emptyset,$$ obtaining the result that $C(K)=\emptyset$.

In our case, we consider the curve $C_D$ together with a map
$\phi:C_D \to E^{(1)}$, where $E^{(1)}$ is the curve given by the
Weierstrass equation $y^2=x(x+2)(x+6)$. The curve $E^{(1)}$ has
Mordell-Weil group $E^{(1)}(\bQ)$ generated by the 2-torsion
points and $P:=(6,24)$.

\begin{lem} Let $D$ be a square-free integer, and consider the
curve $C_D$, together with the map $\phi:C_D\to E^{(1)}$ defined
as
$$\phi([x_0:x_1:x_2:x_3:x_4]):=\left( \frac{6x_0^2}{x_4^2}, \frac
{24x_0x_1x_2}{x_4^3}\right).$$ Let $P:=(6,24)\in E^{(1)}(\bQ)$.
Then
$$\phi(C_D(\bQ))\subset H:=\{kP\ | \ k \mbox{ odd }\}.$$
\end{lem}

\begin{pf} This lemma is an easy application of the 2-descent.
The map $\phi$ is the composition of two maps. First, the
forgetful map from $C_D$ to the genus one curve in $\bP^3$ given
by the equations
$$
\left\{
\begin{array}{l}
F_{014}:=3X_0^2-4X_1^2+X_4^2=0,\\
F_{024}:=X_0^2-2X_2^2+X_4^2=0,
\end{array}
\right.
$$
given by sending $[x_0:x_1:x_2:x_3:x_4]$ to $[x_0:x_1:x_2:x_4]$.
Multiplying $F_{014}$ by 2 and $F_{024}$ by 6 we obtain the equations
of a 2-descendent
$$
\left\{
\begin{array}{l}
 6X_0^2-2(2X_1)^2=-2X_4^2,\\
 6X_0^2-12X_2^2=-6X_4^2.
 \end{array}
\right.
$$
The second map is the corresponding $4$-degree map $\phi_{6,2}$ from
these curves to $E^{(1)}$ given by the equations above, and determines that
the element $(6,2)\in S^{(2)}(E^{(1)}/\bQ)$, so $\phi(C_D(\bQ))$
is contained in the subset of elements $(x,y)$ of $E^{(1)}(\bQ)$
with $\psi((x,y)):=(x,x+2)=(6,2)$ in $\bQ^*/(\bQ^*)^2$. But
$P:=(6,24)\in E^{(1)}(\bQ)$ is a generator of
$E^{(1)}(\bQ)/E^{(1)}(\bQ)[2]$, and has $\psi(6,24)=(6,2)$, hence
any such point $(x,y)$ is an odd multiple of $P$.
\end{pf}

\

For any prime $q$, we will denote by $H_D^{(q)}\subset H$ the
subset corresponding to
$$H_D^{(q)}:=\red_{q}^{-1}\left(\phi_{q}(C_D(\bF_q)) \cap
\red_{q}(H)\right).$$

First, consider the reduction modulo a prime $q$ dividing $D$, so
a prime of bad reduction. Suppose we have a point
$[x_0:x_1:x_2:x_3:x_4]$ of $C_D$, so $x_0^2$, $x_1^2$, $x_2^2$,
$Dx_3^2$ and $x_4^2$ are coprime integers in arithmetic
progression. By reducing modulo $q$ one gets that $x_0^2$,
$x_1^2$, $x_2^2$, $0$ and $x_4^2$ are in arithmetic progression
modulo $q$, so, after dividing by $x_4^2$, we may suppose that it
is the arithmetic progression $-3,-2,-1,0,1$.

\begin{prop}\label{MWSatD}
Let $q>3$ be a prime number dividing $D$. Then $$H_D^{(q)}=\{kP\ | \ k \mbox{ odd and } x(kP)\equiv -18
\pmodd{q} \},$$
and $H_D^{(q)}$ is independent on $D$.
\end{prop}

\begin{pf}This is an easy application of the ideas above.
Since the only points in the reduction of $C_D$ are the ones
having $x_0^2=-3$, $x_1^2=-2$, $x_2^2=-1$ and $x_4^2=1$, the set
$\phi_q(C_D(\bF_q))$ only contains at most the two points having
$x$-coordinate equal to $6(-3)=-18$.
\end{pf}

\begin{cor}\label{cMWSatD} Suppose that $q>3$ is a prime number such that
$\red_q(H)$ contains a point $Q$ with $x(Q)\equiv -18 \pmodd{q}$.
Then infinitely many pairs of square-free integers $D$
and primitive tuples $[x_0:x_1:x_2:x_3:x_4]\in C_{D}(\bQ)$ exist, such
that either $q$ divides $D$ or $x_3\equiv 0 \pmodd{q}$ .
\end{cor}

\begin{pf}  Let $O_q$ be the order
of $P$ modulo $q$, and let $k$ be such that $x(kP)\equiv -18
\pmodd{q}$. Then $x(k'P)\equiv-18 \pmodd{q}$ for all $k'\equiv k
\pmodd{O_q}$. So, if $k$ is odd or $O_q$ is odd, $H^{(q)}$ has
infinitely many elements. For any point $Q\in H^{(q)}$, we have that
$x(Q)=6 z^2$, for certain $z\in\bQ$, such that $z^2\equiv -3
\pmodd{q}$. Write $z=a/b$ with $a$ and $b\in \bZ$ and coprime.
Then, if we denote by $r:=(a^2-b^2)/4$, then $r\in \bZ$ and
$x_i:=a^2+ir$ are squares for $i=0,1,2$ and $4$, and $a^2+3r
\equiv 0 \pmodd{q}$. Define $D$ as the square-free part of $a^2+3r$, and
we obtain the result by defining $x_3$ such that $a^2+3r=D x_3^2$.
\end{pf}

\

Observe, however, that we do not obtain that
$C_q(\bQ)\ne\emptyset$ for the primes satisfying the hypothesis of the previous corollary. For example, the prime
$q=457$ verifies the conditions of the corollary, but we will show
that $C_{457}(\bQ)= \emptyset$.

Now we will consider primes $q>3$ that do not divide $D$, and are hence
good reduction primes. We will obtain conditions depending on $D$
being a square or not modulo $q$.

\begin{prop}  Let $q>3$ be a prime number that does not divide $D$.
Then $H_D^{(q)}\subset E^{(1)}(\bQ)$ depends only on the Legendre
symbol $(D/q)$. If we denote by $H^{(q),(D/q)}$ the subgroup
corresponding to any $(D/q)$, and by $O_q$ the order of $P\in
E^{(1)}(\bQ)$ modulo $q$, we have that subsets
$M_1^{(q)}$ and $M_{-1}^{(q)}$ of $\bZ/O_q\bZ$ exist, such that
$$H^{(q),(D/q)}=\{kP \ | \ k \mbox{ odd and } \exists m\in M_{(D/q)}^{(q)}
\mbox{ such that } k\equiv m \pmodd{O_P}\}.$$ Moreover, $1\in
M_1^{(q)}$ for any $q>3$, and if $k\in M_{(D/q)}^{(q)}$, then
$-k\in M_{(D/q)}^{(q)}$.
\end{prop}

\begin{pf} First we show that $H^{(q)}_D$ only depends on $(D/q)$.
Suppose that $D\equiv D'a^2 \pmodd{q}$, for certain $a\ne 0 \in
\bF_q$. Then the morphism given by
$\theta([x_0:x_1:x_2:x_3:x_4])=[x_0:x_1:x_2:x_3a^2:x_4]$
determines an isomorphism between $C_{D'}$ and $C_{D}$ defined
over $\bF_q$, clearly commuting with $\phi$, which does not
depend on the $x_3$.

In order to define $M_{(D/q)}^{(q)}$, one computes
$\phi_q(C_D(\bF_q))$ and then intersects it with the subset of
$E^{(1)}(\bF_q)$ of the form $\{kP \ | \ k \mbox{ odd }\}$. Then
$$M_{(D/q)}^{(q)}:=\{k\in \bZ/O_q\bZ\ |\ kP\in
\phi_q(C_D(\bF_q))\}.$$

So $k$ belongs to $M_{(D/q)}^{(q)}$ if some
$Q:=[x_0:x_1:x_2:x_3:x_4]\in C_{D}(\bF_q)$ exists, such that $\phi(Q)=kP$.
But then $\phi([-x_0:x_1:x_2:x_3:x_4])=-kP$.

Finally, if $(D/q)=1$, we can suppose $D\equiv 1 \pmodd{q}$. But
then $Q_0:=[1:1:1:1:1]\in C_{D}(\bF_q)$, and $\phi(Q_0)=P$.
\end{pf}

\

The following table shows some examples of $M_{\pm 1}^{(q)}$ for $5<q<30$ prime.
\begin{center}
\begin{tabular}{|c|c|c|c|}
\hline
$q$ & $O_q$ &  $M_1^{(q)}$ & $M_{-1}^{(q)}$\\
\hline
$7$ &  $6$ & $\{\pm 1\}$ & $\{3\}$\\
\hline
$11$ &  $8$ & $\{\pm 1\}$ & $\{\pm 3\}$\\
\hline
$13$ &  $6$ & $\{\pm 1\}$ & $\{3\}$\\
\hline
$17$ &  $6$ & $\{\pm 1,3\}$ & $\{\,\,\}$\\
\hline
$19$ &  $8$ & $\{\pm 1\}$ & $\{ \pm 3\}$\\
\hline
$23$ &  $3$ & $\{1,2,3\}$ & $\{\,\,\}$\\
\hline
$29$ &  $16$ & $\{\pm 1\}$ & $\{\pm 3,\pm 5,\pm 7\}$\\
\hline
\end{tabular}
\end{center}
We are going to use the previous result to obtain conditions on $D$.

\begin{cor}\label{MW}
 If $C_D(\bQ)\ne \emptyset$ then $D$ satisfies the following conditions:
\begin{itemize}
\item[(i)] $D$ is a square modulo 17, 23, 41, 191, 281, 2027, 836477.
\item[(ii)] $(D/7)=(D/13)$, $(D/11)=(D/19)=(D/241)$, $(D/47 )=(D/73)$,
$(D/149)=(D/673)$, $(D/43)=(D/1723)$, $(D/175673)=(D/2953)$,
$(D/97)=(D/5689)=(D/95737)$, $(D/577)=(D/2281)$,
$(D/83)=(D/4391)=(D/27449)$, $(D/67)=(D/136319)$,
$(D/2111)=(D/2521)$.
\item[(iii)] If $(D/29)=1$ then $(D/11)=1$. If
$(D/149)=1$ then $(D/31)=1$. If $(D/7019)=1$ then $(D/8123)=1$. If
$(D/617)=1$ then $(D/37)=1$, and in this case $(D/7)=1$.
\item[(iv)] If
$(D/83)=-1$ then $(D/11)=-1$. If $(D/2347)=-1$ then $(D/47)=-1$.
If $(D/10369)=-1$ then $(D/47)=-1$.
\end{itemize}
\end{cor}

\begin{pf}
We have computed the sets $M_1^{(q)}$ and $M_{-1}^{(q)}$ for $q <
10^6$ and $O_q\leq 200$. Then the algorithm to obtain the
conditions for the statement is as follows: fix an integer $k\leq
200$ and compute the primes $q$ such that $O_q=k$ and $5<q<10^6$.
For these primes compute $M_1^{(q)}$ and  $M_{-1}^{(q)}$. If
$M_{-1}^{(q)}$ is empty, then $(D/q)=1$ and we obtain (i). If
these sets are equal for different primes, then we obtain (ii).
Now, for any integer $m>1$ such that $mk\le 200$, compute the
primes $p<10^6$ such that $O_p=m k$. Compute $M_1^{(p)}$ and
$M_{-1}^{(p)}$. Now check if $M_1^{(p)}$ (resp. $M_{-1}^{(p)}$)
mod $k$ is equal to some of the sets $M_1^{(q)}$ (resp.
$M_{-1}^{(q)}$) computed above. If this occurs, then we obtain the
rest of the conditions.

For example, looking at the previous table we see that
$M_{-1}^{(17)}=\{\}$, therefore $(D/17)=1$. Now, $O_7=O_{13}$,
$M_{1}^{(7)}=M_{1}^{(13)}$ and $M_{-1}^{(7)}=M_{-1}^{(13)}$ so we
have $(D/7)=(D/13)$. Finally, $O_{29}=2 O_{11}$ and $M_{1}^{(29)}$
mod $O_{11}$ is equal to $M_{1}^{(11)}$ and then we obtain that if
$(D/29)=1$ then $(D/11)=1$.
\end{pf}

\section{Computing all the points for $D=409$}\label{sec_409}

We want to find all the rational points of the curve $C_D$ when we
know there are some. We will concentrate at the end in the case
$D=409$, which is the first number that passes all of the tests (see
Corollary \ref{computations}), but for the main part of the
section we can suppose that $D$ is any prime integer fulfilling the
conditions in Proposition \ref{local}. Observe first that we do have
the $16$ rational points $[\pm 7,\pm13,\pm17,1,\pm23]\in
C_{409}(\bQ)$. Our aim is to show that there are no more.

In recent years, some new techniques have been developed in order
to compute all the rational points of a curve of genus greater
than one over $\bQ$. These techniques work only under some special
hypotheses. For example, Chabauty's method (see \cite{Chabauty},
\cite{Coleman}, \cite{Flynn2}, \cite{Stoll2}, \cite{Stoll},
\cite{McCallum:Poonen}) can be used when the Jacobian of the curve
has rank less than the genus of the curve, or even when there is a
quotient abelian variety of the Jacobian with rank less than its
dimension. In our case, however, the Jacobian of the curve $C_{D}$
is isogenous to a product of elliptic curves, each of them with
rank one or higher (in fact, the Jacobian of $C_D$ must have rank
$\ge 8$ by Proposition \ref{rank}). So we cannot apply this
method. Other methods, like the Dem'janenko-Manin's method \cite{Demjanenko,Manin}, cannot be
applied either. We will instead apply the covering collections
technique, as developed by Coombes and Grant \cite{CG}, Wetherell
\cite{We} and others, and specifically a modification of what is
now called the elliptic curve Chabauty method developed by  Flynn and
Wetherell in \cite{FW} and by Bruin in \cite{Br,Br0}.

The idea is as follows: suppose we have a curve $C$ over a number
field $K$ and an unramified map $\chi:C'\to C$ of degree greater than one, and
defined over $K$. We consider the distinct unramified
coverings $\chi^{(s)}:C'^{(s)}\to C$ formed by twists of the given
one, and we obtain that
$$C(K)=\bigcup_{s} \chi^{(s)}(C'^{(s)}(K)),$$
the union being disjoint. In fact, only a finite number of twists
do have rational points, and the finite (larger) set of twists
with points locally everywhere can be explicitly described. Now
one hopes to be able to compute the rational points of all the
curves $C'^{(s)}$, and therefore also of the curve $C$.

We will consider degree 2 coverings of $C_D$. To construct such
coverings, we will use the description given by Bruin and Flynn in
\cite{BF} of the $2$-coverings of curves which are $2$-coverings
of the projective line. In our case, $C_D$ is not of such form,
but a quotient of $C_D$ is of this form. Therefore we will use a
$2$-covering for such a quotient. Specifically, we will use one of
the five genus 1 quotients, particularly the quotient
$$
F^{(4)}_D\,:DX_3^2 =t^4-8t^3+2t^2+8t+1\,,
$$
along with the forgetful map $\phi^{(4)}:C_D\longrightarrow
F_D^{(4)}$ given by
$t= \frac{X_0-X_1}{X_2-X_1}$.

Observe first that the curve $C_D$ has some $\bQ$-defined
automorphisms $\tau_i$ of order 2, defined by sending
$\tau_i(x_j)=x_j$ if $j\ne i$, $\tau_i(x_i)=-x_i$. All them,
together with their compositions, form a subgroup $\Upsilon$ of
the automorphisms isomorphic to $(\bZ/2\bZ)^4$. For every
$\bQ$-defined point of $C_D$, composing with these automorphisms
gives 16 different points. Given $Q\in C_D(\bQ)$, we denote by
$T_Q$ the set of all of these $16$ different points. Observe that
$\phi^{(4)}(T_Q)$ is formed by $8$ distinct points.

\begin{lem}\label{involutions}
The involutions $\tau_0$, $\tau_1$, $\tau_2$ and $\tau_3$ give the following involutions on $F^{(4)}_D$:
$$
\tau_0(t,X_3)=\left(\frac{1-t}{1+t},\frac{2X_3}{(1+t)^2}\right), \,
\tau_1(t,X_3)=\left(\frac{-1}{t},\frac{X_3}{t^2}\right),\,
\tau_2(t,X_3)=\left(\frac{t+1}{t-1},\frac{2X_3}{(t-1)^2}\right), \,
$$
and $\tau_3(t,X_3)=(t,-X_3)$. Moreover, if $F_D^{(4)}(\bQ)\ne
\emptyset$ and $\psi:F_D^{(4)}\rightarrow E_D^{(0)}$ is an
isomorphism, then the involutions of $E_D^{(0)}$ given by
$\epsilon_i:=\psi\tau_i\tau_3\psi^{-1}$ for $i=0,1,2$, are
independent of $\psi$. Specifically, $\epsilon_i=\epsilon_{R_i}$
for $R_0=(0,0)$, $R_1=(-D,0)$ and $R_2=(-4D,0)$, where
$\epsilon_Q$ denotes the translation by $Q\in E_D^{(0)}$.
\end{lem}

\begin{pf}
The computation to check the formulae for the involutions on $F^{(4)}_D$ is straightforward.

First, we show that the involutions $\epsilon_i$ are independent of the fixed
isomorphism $\psi$. In order to show this, recall that, in any
elliptic curve, any involution $\epsilon$ that has no fixed points
should be of the form $\epsilon_R(S)=S+R$, for a fixed $2$-torsion
point $R$. Since $\tau_i\tau_3$ has no fixed points in $F^{(4)}_D$, their
corresponding involution $\epsilon_i$ in $E_D^{(0)}$ must be equal to some
$\epsilon_{R_i}$, hence, determined by the corresponding $2$-torsion point $R_i$,
which is equal to $\epsilon_i(0)$. Now, replacing the isomorphism
$\psi$ from $F^{(4)}_D$ to $E_D^{(0)}$ is equivalent to conjugate $\epsilon_i$ by a
translation $\epsilon_Q$ of $E_D^{(0)}$ with respect to a point $Q$ in
$E_D^{(0)}$, so, in principle we obtain a new involution $\epsilon_{-Q}
\epsilon_i \epsilon_Q$, again without fixed points. But
$\epsilon_{-Q} \epsilon_i \epsilon_Q(0)=\epsilon_{-Q}(\epsilon_i(Q))=\epsilon_{-Q}(Q+R_i)=R_i$, so
$\epsilon_{-Q} \epsilon_i\epsilon_Q=\epsilon_i$.

Second, since $\epsilon_i$ is independent of the chosen
isomorphism $\psi$, and also does not depend on the field $K$, we
can replace with a field $K':=K(\sqrt{D})$ where we have
$F^{(4)}_D\cong F^{(4)}_1$, so we are reduced to the case $D=1$.
In this case, a simple computation by choosing some point in
$F^{(4)}_1(\bQ)$ shows that $\epsilon_i=\epsilon_{R_i}$ where
$R_0=(0,0)$, $R_1=(-1,0)$ and $R_2=(-4,0)$ in $E_1^{(0)}$, which
gives the result when we translate them to the curve $E_D^{(0)}$.
\end{pf}

\

Now, we want to construct some degree two unramified coverings of
$F^{(4)}_D$. All these coverings are, in this case, defined over
$\bQ$, but we are interested in special equations not defined over
$\bQ$. The idea is simple: first, factorize the polynomial
$q(t):=t^4-8t^3+2t^2+8t+1$ as the product of two degree $2$
polynomials (over some quadratic extension $K$). In the sequel of
this section, we will denote $K:=\bQ(\sqrt{2})$. Then we have the
factorization $q(t)=q_1(t)q_2(t)$ over $K$ where
$q_1(t):=t^2-(4+2\sqrt{2})t-3 - 2\sqrt{2}$ and $q_2(t):=
\overline{q_1}(t)$, where $\overline{z}$ denotes the Galois
conjugate of $z\in K$ over $\bQ$. We could have chosen other
factorizations over other quadratic fields, but this one is
especially suitable for our purposes as we will show in the
sequel. Then, for any $\delta\in K$, the curves $F'_{\delta}$
defined in $\bA^3$ by the equations
$$
F'_{\delta} :
\left\{
\begin{array}{rcl}
\delta y_1^2 &=q_1(t)=& t^2-(4+2\sqrt{2})t-3 - 2\sqrt{2}\\
(D/\delta) y_2^2 &=q_2(t)=& t^2-(4-2\sqrt{2})t-3 + 2\sqrt{2}
\end{array}
\right.
$$
along with the map $\nu_{\delta}$ that gives $X_3=y_1y_2$ are
all the twists of an unramified degree two coverings of $F^{(4)}_D$.
Observe that, for any $\delta$ and $\delta'$, such that
$\delta\delta'$ is a square in $K$, we have an isomorphism between
$F'_{\delta}$ and $F'_{\delta'}$. So we only need to consider the
$\delta$'s modulo squares. This also means that we can suppose
that $\delta\in \bZ[\sqrt{2}]$. However, only very few of them are
necessary in order to cover all the rational points of $F^{(4)}_D$. A
method to show this type of results is explained in \cite{BF}, but
we will follow a different approach.

\begin{lem}\label{Delta}
Let $D>3$ be a prime number such that
$F_D^{(4)}(\bQ)\neq\emptyset$. Let $\alpha\in\bZ[\sqrt{2}]$ be such that
$\nu_{\alpha}(F'_{\alpha}(K))\cap F^{(4)}_D(\bQ)\ne \emptyset$, then
$$F^{(4)}_D(\bQ) \subset \nu_{\alpha}(F'_{\alpha}(K))\cup
\nu_{\overline{\alpha}}(F'_{\overline{\alpha}}(K)) \cup
\nu_{-\alpha}(F'_{-\alpha}(K))\cup
\nu_{-\overline{\alpha}}(F'_{-\overline{\alpha}}(K)) .
$$
Moreover, for any $Q\in C_D(\bQ)$, either
$$
\phi^{(4)}(T_Q)\cap \nu_{\alpha}(F'_{\alpha}(K))\ne \emptyset\quad\mbox{or}\quad\phi^{(4)}(T_Q)\cap
\nu_{-\overline{\alpha}}(F'_{-\overline{\alpha}}(K)))\ne
\emptyset.
$$
\end{lem}

\begin{pf}
Observe that, for any point $P\in F^{(4)}_D$, an easy calculation
shows that
$$q_1(t(\tau_0(P)))=\frac{2}{(1+t(P))^2} q_1(t(P)) \mbox{ and }
q_1(t(\tau_1(P)))=-\frac {(1+\sqrt{2})^2}{(t(P))^2}q_2(t(P)),$$
where $t(R)$ denotes the $t$-coordinate of the point $R$. This implies that, if $P$ is in $\nu_{\alpha}(F'_{\alpha}(K))\cap
F^{(4)}_D(\bQ)$, then $\tau_0(P)$ and $\tau_3(P)$ also are, and
$\tau_1(P)$ and $\tau_2(P)$ are in
$\nu_{-\overline{\alpha}}(F'_{-\overline{\alpha}}(K))\cap
F^{(4)}_D(\bQ)$. This last fact shows the last assertion of the
lemma.

Now, using a fixed point $P\in F^{(4)}_D(\bQ)$, we choose
$\alpha\in \bZ[\sqrt{2}]$ such that $P\in
\nu_{\alpha}(F'_{\alpha}(K))$, and an isomorphism $\psi_P$ of
$F^{(4)}_D$ with its Jacobian $E:=E_D^{(0)}$, by sending $P$ to
$0$ (this isomorphism is determined, modulo signs, by this fact).
Via this isomorphism, one can identify the degree two unramified
covering $\nu_{\alpha}$ with a degree two isogeny
$\widetilde{\nu}:E'\to E$. Recall that $E$ can be written by the
Weierstrass equation $y^2=x^3+5Dx^2+4D^2x$, and that the degree
two isogenies are determined by a non-trivial $2$-torsion point.

By Lemma \ref{involutions}, we have
$\psi_P(\tau_0\tau_3(P))=\epsilon(0)=(0,0)$. But $\tau_0\tau_3(P)$
also belongs to $\nu_{\alpha}(F'_{\alpha}(K))$, and hence $(0,0)$
must be in $\widetilde{\nu}(E'(\bQ))$, thus determining the
isogeny as the one corresponding to $(0,0)$.

Now we use the standard descent via a $2$-isogeny. One obtains that
$E(\bQ)/\widetilde{\nu}(E'(\bQ))$ is injected inside the subgroup
of $\bQ^*/(\bQ^*)^2$ generated by $-1$ and the prime divisors of
$4D^2$.  Since $D$ is prime, the only possibilities are $-1$, $2$
and $D$, which become only $-1$ and $D$ over $K^*/(K^*)^2$. Hence,
we need only four twists of $\widetilde{\nu}$ over $K$ in order to
cover all the points of $E(\bQ)$. Note that the twist corresponding to
$1$ is identified with $\nu_{\alpha}$. To find the twist corresponding to $-1$ one can argue in the following way: when replacing the field $K$ with $K(\sqrt{-1})$ then $-1$ becomes equal to $1$ modulo squares and not $D$ or $-D$, and the same applies to $\alpha$ and $-\alpha$. Hence $-1$ is identified to $\nu_{-\alpha}$. A similar argument, but postulating that $\alpha\overline{\alpha}$ is equal to $D$ modulo squares in $K$, shows that $D$ corresponds to $\nu_{\overline{\alpha}}$.
\end{pf}

\

In order to obtain some coverings of $C_D$ from these coverings of $F^{(4)}_D$  we write $C_D$ in a different form, the one given by the
following equations in $\bA^3$:
\begin{equation}
C_D:\,\{\,\,DX_3^2=q(t) \ , \  X_4^2=p(t)\,\,\}
\end{equation}
where $p(t)=t^4
- 12 t^3 + 2 t^2 + 12 t + 1$. Then, the lemma above implies that
any rational point of $C_D$, modulo the automorphisms in
$\Upsilon$, comes from a point in $K$ of one of the curves
$C'_{\delta}$, with $\delta=\alpha$ or $\delta=-\alpha$, given by the
following equations in $\bA^4$:
$$
C'_\delta:\,\{\,\, \delta y_1^2=q_1(t) \ , \  (D/\delta)
y_2^2=q_2(t) \ , \  X_4^2=p(t)\,\,\}$$ (and, moreover, with $t\in
\bQ$)  by the natural map $\mu_\delta$. Observe, before
continuing, that any rational point in $C_D$ comes from a point in
the affine part in the previous form, which is singular at infinity,
since $D$ is not a square in $\bQ$.

Now we consider the following hyperelliptic quotient $H_{\delta}$
of the curve $C'_{\delta}$, which can be described by the equation
$$
H_\delta\,:\,\delta W^2=q_1(t)p(t),
$$ and where the quotient map $\eta$ is determined
by saying that $W=y_1X_4$.

The following lemma is a simple verification.

\begin{lem}
Let $E_{\delta}$ be the elliptic curve defined by the equation
$$
E_{\delta}\,:\,\delta y^2 = x^3+5\sqrt{2}x^2-x.
$$
Then, the following equation determines a non-constant morphism
from the genus $2$ curve $H_{\delta}$ to $E_{\delta}$:
$$
\varphi: H_{\delta} \to E_{\delta}\,,\quad
\varphi(t,W)=\left( \frac {-2(-3+2\sqrt{2}) q_1(t)}{(t-\sqrt{2}+1)^2}, \frac{3(-4 + 3 \sqrt{2}) W} {(t-\sqrt{2}+1)^3}\right).
$$
\end{lem}

\begin{rem} The group of automorphism of the genus $2$ curve $H_{\delta}$ is generated by
a non-hyperelliptic involution $\tau$ and by the hyperelliptic
involution $\omega$. Then, we have that the elliptic curve
$E_{\delta}$ is $H_{\delta}/\langle\tau \rangle$. The other
elliptic quotient $E'_{\delta}$ is obtained by $\tau\omega$; that
is, $E'_{\delta}=H_{\delta}/\langle\tau\omega \rangle$. It is easy
to compute that $E'_{\delta}\,:\,\delta y^2 =
x^3+9\sqrt{2}x^2-81x$. Therefore, $\jac{H_{\delta}}$ is
$\bQ(\sqrt{2})$-isogenous to $E_{\delta}\times E'_{\delta}$.
Moreover,  $E_{1}$ and $E'_{1}$ are $\bQ(\sqrt{2})$-isomorphic
respectively to \verb+384f2+ and \verb+384c2+ in Cremona's tables,
so $E_{\delta}$ and $E'_{\delta}$ are $\delta$-twists of them.
\end{rem}

\begin{rem} The fact that $H_{\delta}$ has such elliptic quotient defined over
$K$ is the main reason we consider these specific $2$-coverings of
$C_D$. If we want to carry out the same arguments with other
$2$-coverings, coming from $2$-coverings of $F_D^{(4)}$ or from
$2$-coverings of other genus $1$ quotients $F^{(i)}_D$, we will
not obtain such a quotient defined over a quadratic extension of
$\bQ$.
\end{rem}

In the following proposition we will determine a finite subset of
$E_{\delta}(K)$ containing the image of the points $Q$ in
$C_{\delta}(K)$ such that $\mu_{\delta}(Q)\in C_D(\bQ)$.

\begin{prop} Let $D>3$ be a prime number such that
$C_D(\bQ)\neq\emptyset$. Consider $P\in C_D(\bQ)$. Then $\tau \in \Upsilon$ exists such that $\tau(P)=\mu_{\delta}(Q)$ for
$\delta=\alpha$ or $\delta=-\alpha$, with $Q\in C'_{\delta}(K)$.
Let $R:=\varphi(\eta(Q))\in E_{\delta}(K)$ be the corresponding
point in $E_{\delta}$. Then
$$R\in \{(x,y)\in E_{\delta}(K)\ | \
\pi(x,y):=\frac{2(-4+2\sqrt{2}-x(1-\sqrt{2}))}{(6-4\sqrt{2}-x)}\in
\bQ \}.$$
\end{prop}

\begin{pf} Part of the lemma is a recollection of what we have proved in lemmas
above. Only the last assertion needs a proof. So, suppose we have
a point $Q\in C'_{\delta}(K)$ such that $\mu_{\delta}(Q)\in
C_D(\bQ)$. Then the $t$-coordinate of $Q$ is in $\bQ$, since
$\mu_{\delta}$ leaves the $t$-coordinate unchanged. This implies
that the $x$-coordinate of $R:= \varphi(\eta(Q))$, that is $\frac
{-2(-3+2\sqrt{2}) q_1(t)}{(-1 + \sqrt{2}-t)^2}$, must come from a
rational number $t$. This again implies that the sum of the $t$-coordinates of the two pre-images of $R$ is a rational number. But
this sum can be expressed in the $x$-coordinate of $R$ as
$\pi(x,y).$
\end{pf}

\

The following diagram illustrates all the curves and morphisms
involved in our problem:
$$
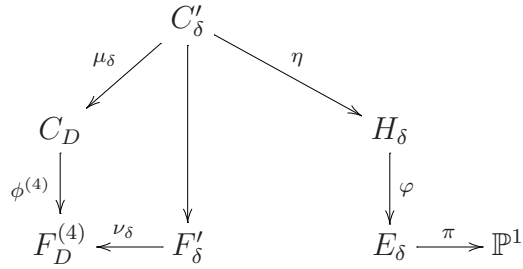
\xymatrix{
    &              C'_\delta \ar@{->}[dd]    \ar@{->}[dl]_{\mu_\delta} \ar@{->}[drr]^\eta   &      &  &  \\
 C_D\ar@{->}[d]_{\phi^{(4)}}  &     &  & H_\delta \ar@{->}[d]^\varphi  & \\
 F_D^{(4)} &F'_\delta  \ar@{->}[l]_{\nu_\delta}& & E_\delta \ar@{->}[r]^\pi  & \mathbb{P}^1
}
$$

Hence, to find all the points in $C_D(\bQ)$, it is enough to find
all the points $(x,y)$ in $E_{\delta}(K)$ such that $\pi(x,y)\in
\bQ$ for $\delta=\alpha$ or $\delta=-\alpha$. But this is what the
so-called elliptic curve Chabauty does, if the rank of the group of
points $E_{\delta}(K)$ is less than or equal to 1. And this seems
to be our case in the cases we are considering.

\begin{exmp}\label{example}
Let us consider the case $D=409$. The 16 points $[\pm
7,\pm13,\pm17,1$, $\pm23]$ give the 8 points in $F^{(4)}_{409}$
with $t\in \{-3/2, -5, 2/3, 1/5 \}$. Take $\alpha:=21+4\sqrt{2}$, which satisfies the hypothesis of Lemma \ref{Delta}. Then the 8 points in $C_{409}$
with $t=-3/2$ and $t=-5$ come from the 16 points in $C_{\alpha}'$
given by $[t,y_1,y_2,X_4]=[-3/2,\pm 1/2,\pm 1/2,\pm 23/4]$ and
$[-5,\pm \sqrt{2},\pm \sqrt{2},\pm 46]$ respectively, which in
turn give the 4 points in $H_{\alpha}$ given by $[t,W]=[-3/2,\pm
23/8]$ and $[-5,\pm 46\sqrt{2}]$. Finally, these $4$ points give the
following 2 points 
$E_{\alpha}$:
$$\left( \frac {-2}{49}(-663+458\sqrt{2}), \pm \frac {69}{343}(-232+163\sqrt{2}) \right).$$
The other points with $t=2/3$ and $1/5$ rise to points in
$E_{-\overline{\alpha}}(K)$, as shown in Lemma \ref{Delta}. We
will show that these points in $E_{\alpha}(K)$ are the only points
$R$ with $\pi(R)\in \bQ$, and that there are no such points in
$E_{-\alpha}(K)$.
\end{exmp}

\subsection{Elliptic curve Chabauty}

In order to apply the elliptic curve Chabauty technique \cite{Br,Br0}, we first need to fix a rational prime $p$ such that it is inert over $K$ and
$E_{\delta}$ has good reduction over such $p$. The smallest such
prime under our conditions is $p=5$, since by Proposition \ref{local} we have $D\equiv \pm 1
\pmodd{5}$. Denote by  $\widetilde{E_{\delta}}$  the reduction
modulo $5$ of  $E_{\delta}$, which is an elliptic curve over
$\bF_{25}:=\bF_5(\sqrt{2})$. Then the elliptic curve Chabauty method
will allow us to bound, for each point $\widetilde{R}$ in
$\widetilde{E_{\delta}}(\bF_{25})$, the number of points $R$ in
$E_{\delta}(K)$ reducing to that point $\widetilde{R}$, and such
that $\pi(R)\in\bQ$, if the rank of the group of points
$E_{\delta}(K)$ is less than or equal to 1. In the next lemma we
will show that, in fact, we only need to consider four (or two)
points in $\widetilde{E_{\delta}}(\bF_{25})$, instead of all $32$ points.

\begin{lem} Let $D$ be a square-free integer such that $D\equiv \pm 1 \pmodd{5}$, and
let $\delta\in \bZ[\sqrt{2}]$ and $Q\in C'_{\delta}(K)$ be such
that $\mu_{\delta}(Q)\in C_D(\bQ)$. Let $R:=\varphi(\eta(Q))\in
E_{\delta}(K)$ be the corresponding point in $E_{\delta}$. Then
$\pi(R) \equiv -1 \pmodd{5}$ or $\pi(R)\equiv \infty \pmodd{5}$.

Moreover, if the rank of the group of points $E_{\delta}(K)$ is
 equal to 1, the torsion subgroup has order $2$, and the reduction of the generator has
order $4$, then only one of the two cases may occur.
\end{lem}

\begin{pf} We repeat the whole construction of the coverings, but
modulo $5$. First, observe that, since $D\equiv \pm 1 \pmodd{5}$,
the only $\bF_5$-rational points of $\widetilde{C_D}$ are the ones
with coordinates $[\pm 1:\pm 1:\pm 1:1:\pm 1]$. So the
$t$-coordinates of these points are $t=0,1,4$ and $\infty$.
Substituting these values in $q_1(t)$ modulo $5$, we always obtain
squares in $\bF_{25}$. This implies that all the twists of the
curves involved are all isomorphic modulo $5$ to the curves with
$\delta=1$.

Consider the curve $\widetilde{H_1}$ over $\bF_{25}$. A simple
computation shows that the only points in $\widetilde{H_1}$ whose
$t$-coordinates are $\bF_{5}$-rational are the points with $t=0$, $t=1$ and
the two points at infinity. Now, the image by
$\varphi$ in $\widetilde{E_1}$ of these points are  equal to the points with
$x$-coordinate equal to $-\overline{\xi}=-1+\sqrt{2}$ in the first
two cases, and equal to $\xi=1+\sqrt{2}$ for the points at
infinity. In the first case we have that $\pi(-1+\sqrt{2})\equiv
-1 \pmodd{5}$, and in the second one we have that
$\pi(1+\sqrt{2})\equiv \infty \pmodd{5}$.

Now, the curve $\widetilde{E_1}$, given by the equation $y^2= x^3
+ 4x$, has $32$ rational points over $\bF_{25}$, and
$\widetilde{E_1}(\bF_{25})\cong \bZ/4\bZ\oplus\bZ/8\bZ$ as abelian
group, with generators some points $P_4$ and $P_8$ with
$x$-coordinate equal to $\xi=1+\sqrt{2}$ and
$\sqrt{2}\xi=2+\sqrt{2}$ respectively. We then obtain that
$$\{R\in\widetilde{E_1}(\bF_{25})\ | \pi(R)=\infty\}=\{P_4,-P_4\}$$
and
$$\{R\in\widetilde{E_1}(\bF_{25})\ | \pi(R)=-1\}=\{2P_8+P_4,-2P_8-P_4\}.$$
Now, if the rank of the group of points $E_{\delta}(K)$ is less
than or equal to 1, the torsion subgroup has order $2$, and the
reduction of the generator has order $4$, then the reduction of
$E_{\delta}(K)$ is a subgroup of $\widetilde{E_1}(\bF_{25})$
isomorphic to $\bZ/4\bZ\oplus\bZ/2\bZ$. But the subgroup generated
by $P_4$ and $2P_8+P_4$ is isomorphic to $\bZ/4\bZ\oplus\bZ/4\bZ$,
and the reduction cannot therefore contain both points.
\end{pf}

\

In order to use elliptic curve Chabauty, it is advisable to
transform the equation that gives $E_{\delta}$ into a Weierstrass
equation, by doing the standard transformation sending $(x,y)$ to
$(\delta x,\delta y)$. We obtain the equation
$$y^2=x^3+ 5\sqrt{2}\delta x^2-\delta^2 x.$$
We will denote this elliptic curve, by abuse of notation, by
$E_{\delta}$. Moreover, the map $\pi$ becomes the map
$f:E_{\delta}\to \bP^1$, given by
$$f(x):=\frac{(2\sqrt{2}-2)x+\delta(4\sqrt{2} -
8)}{\delta(-4\sqrt{2}+6)-x}.$$

Let us explain first the idea of the elliptic curve Chabauty method. For
a given $D$, we fix a $\delta=\alpha$ or $\delta=-\alpha$, and we want to
compute the set
$$\Omega_{\delta}:=\{Q\in E_{\delta}(K)\ | \ f(Q)\in \bQ \mbox{
and } f(Q)\equiv -1, \infty \pmodd{5}\}.$$ As we have already remarked,
we need first to compute the rank of the group $E_{\delta}(K)$,
which should be less than or equal to one. We will also need to explicitly know the
torsion subgroup of that group, and some non-torsion point if the
rank is $1$, which is not an $\ell$-multiple of a $K$-rational
point for some primes $\ell$ to be determined (in our cases, they
will only be $\ell=2$). In the cases where we already know some points
in $E_{\delta}(K)$, those coming from the known points in
$C_D(\bQ)$, we will show that those points are non-torsion points.

We have two cases we would like to consider. The first case is when we
will not know any point $R\in E_{\delta}(K)$ such that $f(R)\in
\bQ$. In these cases we hope to show that
$\Omega_{\delta}=\emptyset$ by just proving that the reduction of
the group $E_{\delta}(K)$ does not contain any point
$\widetilde{Q}$ such that $\widetilde{f}(\widetilde{Q})\in \bF_5$.
We do so for the two cases in the following lemma.

\begin{lem}\label{Omegaxi} Take $D=409$ and $\alpha=21+4\sqrt{2}$. Then the
elliptic curves $E_{\alpha}$ and $E_{-\alpha}$ have rank 1 over
$K$ and torsion subgroup isomorphic to $\bZ/2\bZ$ (generated by the point $(0,0)$). The points $P=((-30\sqrt{2}- 43)/2
,(759\sqrt{2}+ 1104)/4)$ in $E_{\alpha}(K)$ and the point $P'\in
E_{-\alpha}(K)$ with $x$-coordinate equal to
$$\frac{29769295809708\sqrt{2} + 42339835565318}{4185701809},
$$
generate the free part of the corresponding Mordell-Weil group.

Moreover, if $R\in \Omega_{-\alpha}$ then $f(R)\equiv  \infty \pmodd{5}$ and if $R\in \Omega_{\alpha}$ then $f(R)\equiv -1 \pmodd{5}$.
\end{lem}

\begin{pf}
The first part of the previous statement was obtained by using the
\verb|MAGMA| function \verb|DescentInformation|. For our elliptic
curves  $E_{\alpha}$ and $E_{-\alpha}$, this function has
unconditionally computed that the rank of both elliptic curves is
$1$, and it has returned the generators of these  Mordell-Weil
groups.

The last assertions are shown by proving that the subgroup
generated by the reduction modulo $5$ of the point $P'$ and the
point $(0,0)$ does not contain any point with image by
$\widetilde{f}$ equal to $-1$, and that the subgroup generated by
the reduction modulo $5$ of the point $P$ and the point $(0,0)$
does not contain any point with image by $\widetilde{f}$ equal to
$\infty$. These last two cases are in fact instances of the
previous lemma, since the reduction of the points $P$ and $P'$
have order $4$.
\end{pf}

\

Now, in order to show that $\Omega_{-\alpha}$ is, in fact, empty, we
need to use information from some other primes. This is what we do in
the following lemma.

\begin{lem}\label{Omegaminus} Take $D=409$ and $\alpha=21+4\sqrt{2}$. Then
$\Omega_{-\alpha}=\emptyset$.
\end{lem}

\begin{pf} By using reduction modulo $5$, we obtain that any point $R$ in
$\Omega_{-\alpha}$ should be of the form $R=(4n+1)P'+(0,0)$ for some
$n\in \bZ$, since it should reduce to the point $\widetilde{P'+T}$,
and the order of $\widetilde{P'}$ is 4.

Now we reduce modulo $13$. One shows easily that the order of $P'$
modulo $13$ is equal to $24$, and that the points $R\in
E_{-\alpha}(K)$ such that $f(R)\in \bP^1(\bQ)$ reduce to the points
$6P'$ or $12P'+(0,0)$. Hence the points $R$ should be of the form
$R=(24n+6)P'$ or $(24n+12)P'+(0,0)$. Comparing with the result
obtained from the reduction modulo 5, we obtain the result that there is no such
point.
\end{pf}

\

The second case is where we already know some points $R\in
\Omega_{\delta}$. Then our objective will be to show there are no
more, by showing that the set $$\Omega_{\delta,R}:=\{Q\in
E_{\delta}(K)\ | \ Q\in \Omega_{\delta} \mbox{ and } Q\equiv R
\pmodd{5}\}$$ only contains the point $R$. This is done by
translating the problem of computing the number of points in
$\Omega_{\delta,R}$ into a problem of computing the number of
$p$-adic zeros of some formal power series, and using Strassmann's
Theorem to do so.

\begin{prop}
Take $\alpha=21+4\sqrt{2}$, and consider the point
$$
R=\left(
\frac {-2}{49}(-663+458\sqrt{2})\alpha, \frac{69}{343}(-232+163\sqrt{2}) \alpha^2
\right).
$$
Then
$$\Omega_{\alpha}=\{Q\in E_{\alpha}(K)\ | \ f(Q)\in \bQ \mbox{
and } f(Q)\equiv -1 \pmodd{5}\}=\{R,-R\}.$$
\end{prop}

\begin{pf}
First of all, observe that the order of the reduction of $P$ modulo
$5$ is $4$. Also, any point $R'$ in $\Omega_{\alpha}$ reduces modulo
5 to the points $\pm R$, so it is of the form $\pm R+4nP$. We are
going to prove there is only one point in $\Omega_{\alpha}$ reducing
to $R$, and we deduce the other case by using the $-1$-involution.

Observe that any point in $E_{\alpha}(K)$ that reduces to $0$ modulo
$5$ is of the form $4nP$ for some $n\in \bZ$. We are going to
compute the $z$-coordinate of that points, where $z=-x/y$ if
$P=(x,y)$, as a formal power series in $n$. Denote by $z_0$ the
$z$-coordinate of $4P$. The idea is to use the formal logarithm
$\log_{E}$ and the formal exponential $\exp_{E}$ of the formal group
law associated to $E_{\alpha}$. These are formal power series in
$z$, one inverse to the other insofar as the composition is concerned, and
such that
$$
\log_E(z\mbox{-coord}(G+G'))=\log_E(z\mbox{-coord}(G))+
\log_E(z\mbox{-coord}(G'))
$$
for any $G$ and $G'$ reducing to $0$ modulo $5$, and where the power
series are evaluated in the completion of $K$ at $5$. Thus, we obtain
that
$$
z\mbox{-coord}(n(4P))=\exp_{E}(n\log_{E}(z_0)),
$$
which is a power series in $n$.

Now, we are going to compute $f(R+4nP)$ as a power series in $n$. To
do so, we use that, by the addition formulae,
$$
x\mbox{-coord}(R+G)=\frac
{w(z)(1+y_0w(z))^2-(a_2w(z)+z+x_0w(z))(z-x_0w(z))^2}{w(z)(z-x_0w(z))^2}
$$
where $R=(x_0,y_0)$, $a_2=5\sqrt{2}\alpha$, $z$ is the
$z$-coordinate of a point $G$ reducing to $0$ modulo $5$ and
$w(z)=-1/y$ evaluated as a power series in $z$. This function is a
power series in $z$, starting as
$x\mbox{-coord}(R+G)=x_0+2y_0z+(3x_0^2+2a_2x_0+a_4)z^2+O(z^3)$,
where $a_4=-\alpha^2=y^2/x-(x^2+ 5\sqrt{2}\alpha x)$. Hence we
obtain that $f(R+4nP)=f(x\mbox{-coord}(R+n(4P))$ can be expressed
as a power series $\Theta(n)$ in $n$ with coefficients in $K$. We
express this power series as
$\Theta(n)=\Theta_1(n)+\sqrt{2}\Theta_2(n)$, with $\Theta_i(n)$
now being a power series in $\bQ$. Then $f(R+4nP)\in \bQ$ for some
$n\in \bZ$ if and only if $\Theta_2(n)=0$ for that $n$. Observe
also that, since $f(R)\in \bQ$, we will obtain that
$\Theta_2(0)=0$, so $\Theta_2(n)=j_1n+j_2n^2+j_3n^3+\cdots$. To
conclude, we will use Strassmann's Theorem: if the $5$-adic
valuation of $j_1$ is strictly smaller than the $5$-adic valuation
of $j_i$ for any $i>1$, then this power series has only one zero
at $\bZ_5$, and this zero is $n=0$. In fact, one can easily show
that this power series verifies that the $5$-adic valuation of
$j_i$ is always greater or equal to $i$, so, if we show that
$j_1\not{\!\!\equiv} 0 \pmodd{5^2}$ we have concluded.

In order to do all this explicitly, we will work modulo some power
of $5$. In fact, working modulo $5^2$ will be sufficient. We have
that $z_0=z\mbox{-coord}(4P)\equiv -10\sqrt{2} + 5 \pmodd{5^2}$,
and that $z\mbox{-coord}(n(4P))\equiv (15\sqrt{2} + 5)n
\pmodd{5^2}$. Finally, we obtain that $\Theta(n)\equiv
19+(15\sqrt{2}+20)n \pmodd{5^2}$, hence $\Theta_2(n)\equiv 15n
\pmodd{5^2}$, so $j_1\equiv 15 \pmodd{5^2}$ which completes the
proof.
\end{pf}

\

An alternative way of proving this result is to use the built-in
\verb+MAGMA+ function \verb+Chabauty+.
The answer is that there are only 2 points $R'$ in $E_{\alpha}(K)$
such that $f(R')\in \bQ$, both having $f(R')=13/2$. Since we
already have two points $\pm R$, both giving $f(R)=13/2$, we have
concluded.

\section{Explicit computations and conjectures}\label{sec9}

We have followed two different approaches to compute for which
square-free integers $D$ there are non-constant arithmetic
progressions of five squares over $\bQ(\sqrt{D})$. On one hand,
for each $D$ we have checked if $D$ passes all the sieves from the
previous sections, obtaining the following result.

\begin{cor}\label{computations}
Let $D<10^{13}$ be a square-free integer such that $C_D(\bQ)\ne
\emptyset$, then $D=409$ or $D=4688329$.
\end{cor}

\begin{pf}
First, for each $D$ we have passed all the local conditions
(Proposition \ref{local}) and the conditions coming from the
Mordell-Weil sieve (Corollary \ref{MW}). Only $1048$ values of $D$
have passed these sieves. To discard all the values except $D=409$
and $D=4688329$, we first apply a test derived from Proposition
\ref{MWSatD}. We test if, for any prime $q$ dividing such $D$,
there is an odd multiple $kP$ of the point $P:=(6,24)\in
E^{(1)}(\bQ)$ reducing to a point with $x$-coordinate equal to
$-18$ modulo q. To explicitly verify this condition, we first compute if there is a point $Q$ in $E^{(1)}(\bF_q)$  with $x$-coordinate equal to
$-18$, the order $O_q$
of $P$ in $E^{(1)}(\bF_q)$ and the discrete logarithm $\log(Q,P)$,
i.e. the number $k$ such that $Q=k P$, if it exists. In case there
is no such $Q$ or there is no such logarithm, then $D$ does not
pass the test. Also in the case that $k$ and $O_q$ are even. In the case that it
passes this first test, we combine this information with the
information from the computation of the $M^{(q)}_D$ for the first
100 primes to discard some other cases.

After this last test there are $34$ values of $D$ that survive,
and we then pass a test based on the ternary forms criterion given
by Proposition \ref{ternary}, by using a short program in
\verb+SAGE+ performed by Gonzalo Tornaria. We check that for these
values $r(D,3x^2+9y^2+16z^2)\ne r(D,x^2+3y^2+144z^2)$. Hence for
those values of $D$, $L(E^{(2)}_D,1)\ne 0$, so the analytic rank
of $E^{(2)}_D$ is zero, hence their rank is also 0.

Only $D=409$ and $D=4688329$ survive all these tests, but for these
values we do have points in $C_D(\bQ)$.
\end{pf}

\

On the other hand, remember that if we take $t= \frac{X_0-X_1}{X_2-X_1}$ then an affine model of $C_D$ is defined by:
$$
C_D:\,\{\,\,DX_3^2=t^4-8t^3+2t^2+8t+1 \ , \  X_4^2=t^4 - 12 t^3 + 2 t^2 + 12 t + 1\,\,\}.
$$
Therefore the curve $F^{(3)}$ that consists of removing the variable $X_3$ on $C_D$ has the equation $F^{(3)}\,:\, X_4^2=t^4 - 12 t^3 + 2 t^2 + 12 t + 1$, and a Weierstrass equation is given by $E^{(1)}:y^2 = x(x+2)(x+6)$. Then we have an isomorphism $\psi:E^{(1)}\longrightarrow F^{(3)}$ defined by
$$
\psi(P)= \left(\frac{6-x}{6+3x-y},\frac{-72-108x-18x^2+x^3+48y}{(6+3x-y)^2}\right),
$$

if $P=(x,y)\neq (-2,0),(-3,- 3),(6,24)$ and $\psi(6,24)=\left(\frac{2}{3},\frac{23}{9}\right)$, $\psi(-2,0)=\infty_1$ and $\psi(-3,-3)=\infty_2$,  where $\infty_{1}$ and $\infty_{2}$ denote the two branches at infinity at the desingularization of $F^{(3)}$ at the unique singular point $[0:1:0]\in\bP^2$. This construction  allows us to construct all the non-constant
arithmetic progressions of five squares over all quadratic fields.
Let $P=(2,-8)$, a generator of the free part of $E^{(1)}(\bQ)$,
and let $n$ be a positive integer. Let $(t_n,z_n)=\psi([n]P)$. Now, consider the square-free
factorization of the number
$$
 t_n^4-8t_n^3+2t_n^2+8t_n+1=D_n w_n^2,
$$
where $D_n\in\bZ$ is square-free, $w_n\in \bQ$. Therefore the
following sequence defines a non-constant arithmetic progression of
$5$ squares over $\bQ(\sqrt{D_n})$:
$$
(-t_n^2-2t_n+1)^2\,\,,\,\,(t_n^2+1)^2\,\,,\,\,(t_n^2-2t_n-1)^2\,\,,\,\,D_n\,
w_n^2\,\,,\,\,z_n^2\,,
$$
and we have points
$Q_n:=[-t_n^2-2t_n+1\,:\,t_n^2+1\,:\,t_n^2-2t_n-1\,:\,w_n\,:\,z_n]\in
C_{D_n}(\bQ)$.

\begin{rem}\label{inf} Observe that the pairs $(D_n, Q_n)$ constructed in
this way are different for different $n$. On the other hand, we
cannot be sure that all the fields $\bQ(\sqrt{D_n})$ are different.
However, we do have an infinite number of integers $D$ such that
$C_D(\bQ)\ne \emptyset$. This is because for any integer $D$, the
curve $C_D$, being of genus $5$ (greater than $1$), always has a finite number of
rational points. Since we do have an infinite number of pairs
$(D_n, Q_n)$ with $Q_n\in C_{D_n}(\bQ)$, we do have an infinite
number of different $D_n$.
\end{rem}

\begin{rem}
If we replace $[n]P$ by $Q\in
\{[n_1]T_1+[n_2]T_2+[m]P_0\,|\,n_1,n_2\in\{0,1\},\,m\in\{n,-n-1\}\}$,
where $T_1=(-2,0)$ and $T_2=(-6,0)$ is a basis of
$E^{(1)}(\bQ)_{\mbox{\tiny tors}}$, we obtain the same arithmetic
progression (up to equivalence). Note that if $n=0$, then we
obtain $D_0=1$ and the previous sequence is the constant arithmetic
progression.
\end{rem}

We summarize in the tables \ref{tab:Dn} and \ref{tab:X0} the
computations that we have carried out using the previous algorithm.
We have normalized the elements of the arithmetic progressions to
obtain integers and without squares in common. We have separated
into two tables. In table  \ref{tab:Dn}, $n$ and the factorization
of $D_n$ appears. In the table \ref{tab:X0}, for each value of
$n$, the corresponding factorization of $X_0$ appear. For all the
values of $n$ computed, we have obtained that the fourth element
of the arithmetic progression is $\sqrt{D_n}$ (in the previous notation, $w_n=1$). That is, if we denote by $r=(D_n-X_0^2)/3$, then
the sequence $\{X_k^2=X_0^2+k\,r \,|\,k\in\{0,\dots,4\}\}$ defines
an arithmetic progression over $\bQ(\sqrt{D_n})$.
\begin{table}[h]
\centering
{\tiny
\begin{tabular}{|c|c|}
\hline
$n$ & $D_n$ \\
\hline
$ 1$ & $409 $   \\
\hline
$2  $ & $4688329  $ \\
\hline
$ 3 $ & $ 457\cdot 548240447113 $ \\
\hline
$ 4 $ & $199554894091303668073201  $   \\
\hline
$  5$ & $4343602906873\cdot 53313950039984189254513  $   \\
\hline
$6  $ & $2593\cdot 9697\cdot 4100179090153\cdot 293318691741678881166926936593$  \\
\hline
$ 7 $ & $330823513952828243573122480536077533156064000139119724642295861921  $   \\
\hline
$8$ & $ 24697\cdot 303049\cdot 921429638596379458921\cdot 291824110407387399760153\cdot 3462757049033071137768291886369 $  \\
\hline
\end{tabular}
}
\vspace{1mm}
\caption{}
\label{tab:Dn}
\end{table}

\begin{table}[h]
\centering
{\tiny
\begin{tabular}{|c|c|c|}
\hline
$n$ & $X_0$ \\
\hline
$ 1$  & $ 7$  \\
\hline
$2  $  & $  47\cdot 89$\\
\hline
$ 3 $  & $31\cdot 113\cdot 577$\\
\hline
$ 4 $  & $7\cdot 176201\cdot 515087  $  \\
\hline
$  5$  & $ 2111\cdot 133967\cdot 1134755801$  \\
\hline
$6  $  & $ 119183\cdot 12622601\cdot 2189366343649 $  \\
\hline
$ 7 $  & $2^{10}\cdot 3\cdot 17\cdot 73\cdot 103787\cdot 112261\cdot 963877\cdot 20581582583  $  \\
\hline
$8$  &  $2^{38}\cdot 3^2\cdot 5\cdot 7\cdot 23\cdot 102179447\cdot 1017098920090613939  $  \\
\hline
\end{tabular}
}
\vspace{1mm}
\caption{}
\label{tab:X0}
\end{table}

One can observe that the size of the $D_n$ we encounter grows very quickly, but we do not know if the $D_n$ constructed in this way
always verify that $D_n<D_{n+1}$. We guess that this condition holds. Even more, the previous table and the Corollary \ref{computations} suggest that, in fact,  there no square-free integer $D$ such that $C_D(\bQ)\ne \emptyset$ and $D_n<D<D_{n+1}$ exists.

If we only use the results in section \ref{sec_local} (Proposition \ref{local})
and section \ref{sec_MW} (Corollary \ref{MW}), we obtain the result that the number of
square-free integers $D$ that pass both tests have positive (but
small) density. This is possibly true if we also use the condition
of the rank, for example Proposition \ref{ternary}, since the
number of twists with positive rank of a fixed elliptic curve
should also have positive density. But we suspect that the number
of actual square-free integers $D$ such that $C_D$ has rational
points should have zero density.

\

{\bf Data:} All the \verb+MAGMA+ and \verb+SAGE+ sources are
available on the first author's webpage.


\


\end{document}